%% file: main.tex
\documentclass[11pt]{article}

\usepackage{fullpage}
\usepackage{authblk}
\usepackage{natbib}
\usepackage{algorithm}
\usepackage{algpseudocode}
\usepackage{amsmath,amsthm,amsfonts}
\usepackage{amsmath}
\usepackage[subnum]{cases}
\usepackage{tcolorbox}
\usepackage{mathrsfs}
\usepackage{amssymb}
\usepackage{color}
\usepackage{mathrsfs}
\usepackage{empheq}
\usepackage{enumitem}
\usepackage{bm}
\usepackage{multirow}
\usepackage{booktabs}
\usepackage{makecell}
\usepackage{graphicx}
\usepackage{subcaption}
\usepackage{comment}
\usepackage{cases}
\usepackage{appendix}
\usepackage{tikz}
\usetikzlibrary{arrows,shapes}
\usepackage[colorlinks= true, linkcolor=red, citecolor=blue, urlcolor=black]{hyperref}
\usepackage{cleveref}
\usepackage{marginnote}
\usepackage{diagbox} 
\usepackage{pgfgantt}

\numberwithin{equation}{section}

\theoremstyle{definition}

\newtheorem{theorem}{Theorem}[section]
\newtheorem{lemma}[theorem]{Lemma}

\newtheorem{defn}[theorem]{Definition}

\DeclareMathOperator*{\argmax}{argmax}
\DeclareMathOperator*{\argmin}{argmin}

\usepackage{booktabs}
\usepackage{tikz}
\usepackage{pgfplots}
\usetikzlibrary{calc,intersections}

\usetikzlibrary{shapes.geometric, arrows, positioning}

\tikzset{
	mybox/.style  = {draw, rectangle, minimum width=4cm, minimum height=0.8cm, text centered, text width=4.4cm,   
		font=\normalsize},
	box/.style  = {draw, rectangle, minimum width=2.0cm, minimum height=0.6cm, text centered, text width=3.0cm,   
		font=\normalsize},
	myarrow/.style = {line width=0.2pt, draw=black, -triangle 60, postaction={draw, line width=0.2pt, shorten >=10pt,-}}
}

\tikzstyle{arrow} = [->, >=stealth, -triangle 60]

\allowdisplaybreaks

\makeatletter
\newcommand{\leqnomode}{\tagsleft@true}
\newcommand{\reqnomode}{\tagsleft@false}
\makeatother

\begin{document}

\title{Understanding the PDHG Algorithm via High-Resolution Differential Equations}
\author[1,2]{Bowen Li}
\author[1,2]{Bin Shi\thanks{Corresponding author: \url{shibin@lsec.cc.ac.cn} } }
\affil[1]{Academy of Mathematics and Systems Science, Chinese Academy of Sciences, Beijing 100190, China}
\affil[2]{School of Mathematical Sciences, University of Chinese Academy of Sciences, Beijing 100049, China}
\date\today

\maketitle

\begin{abstract}
The~\textit{least absolute shrinkage and selection operator} (\texttt{Lasso}) is widely recognized across various fields of mathematics and engineering. Its variant,  the generalized~\texttt{Lasso}, finds extensive application in the fields of statistics, machine learning, image science, and related areas. Among the optimization techniques used to tackle this issue, saddle-point methods stand out, with the \textit{primal-dual hybrid gradient} (\texttt{PDHG}) algorithm emerging as a particularly popular choice. However, the iterative behavior of~\texttt{PDHG} remains poorly understood. In this paper, we employ dimensional analysis to derive a system of high-resolution ordinary differential equations (ODEs) tailored for~\texttt{PDHG}. This system effectively captures a key feature of~\texttt{PDHG}, the coupled $x$-correction and $y$-correction, distinguishing it from the proximal Arrow-Hurwicz algorithm.  The small but essential perturbation ensures that~\texttt{PDHG} consistently converges, bypassing the periodic behavior observed in the proximal Arrow-Hurwicz algorithm. Through Lyapunov analysis, We investigate the convergence behavior of the system of high-resolution ODEs and extend our insights to the discrete~\texttt{PDHG} algorithm. Our analysis indicates that numerical errors resulting from the implicit scheme serve as a crucial factor affecting the convergence rate and monotonicity of~\texttt{PDHG},  showcasing a noteworthy pattern also observed for the~\textit{Alternating Direction Method of Multipliers}  (\texttt{ADMM}), as identified in~\citep{li2024understanding}. In addition, we further discover that when one component of the objective function is strongly convex, the iterative average of~\texttt{PDHG} converges strongly at a rate $O(1/N)$, where $N$ is the number of iterations.
\end{abstract}

%

\input{01_intro}

\input{02_prelim}

\input{03_ode}

\input{04_algorithm}

\input{05_extension}

\input{06_conclu}

\section*{Acknowledgements}
Bowen Li was partially supported by the Hua Loo-Keng scholarship of CAS.  Bin Shi was partially supported by Grant No.YSBR-034 of CAS.

\bibliographystyle{abbrvnat}
\bibliography{sigproc}
\end{document}

%% file: 01_intro.tex
\section{Introduction}
\label{sec: intro}

At the turn of this century, a significant discovery was made regarding the observation that sparse representation in a suitable basis or dictionary can effectively model many real-world signals~\citep{candes2006robust, donoho2006compressed}.  This sparse representation, used for signal reconstruction, requires significantly fewer samples and offers enhanced robustness. The $\ell_1$-regularization, acting as a norm that induces sparsity, began to gain recognition in history. A more general formulation,  which also accounts for noise in the observed signal $b$,  emerged as the~\textit{least absolute shrinkage and selection operator} (\texttt{Lasso})~\citep{tibshirani1996regression}. Given $A \in \mathbb{R}^{m \times d_1}$ representing an $m \times d_1$ matrix and $b \in \mathbb{R}^m$ a vector of  $m$ dimensions,  this formulation is expressed as:
\[
\min_{x \in \mathbb{R}^d} \Phi(x) : = \frac{1}{2}\|Ax - b\|^2 + \lambda\|x\|_1,
\]
where the regularization parameter $\lambda > 0$ serves as a tradeoff between fidelity to the measurements and sensitivity to noise.\footnote{Throughout this paper, the notation $\|\cdot\|$ specifically refers to the $\ell_2$-norm or Euclidean norm, $\|\cdot\|_2$, and the $\ell_1$ norm  is defined as
\[ \|x\|_1 = \sum_{i=1}^{d} |x_i|,\]
for any $x \in \mathbb{R}^d$. It is worth noting that the subscript $2$ is often omitted unless otherwise noted.} However, for a wider range of applications,  the
commonly used form is the generalized \texttt{Lasso}~\citep{tibshirani2011solution}, which can be written as
\[
\min_{x \in \mathbb{R}^{d_1}} \Phi(x) : = \frac{1}{2}\|Ax - b\|^2 + \lambda \|Fx\|_1,
\]
with $F \in \mathbb{R}^{d_2 \times d_1}$ incorporating the total variation regularizer  in lieu of the $\ell_1$ norm.  This generalization significantly expands its application across various domains in machine learning and imagine science~\citep{tibshirani2011solution, chambolle2016introduction}, particularly supporting techniques such as~\textit{total-variation denoising}~\citep{rudin1992nonlinear}, \textit{fused lasso}~\citep{tibshirani2005sparsity}, $\ell_1$ \textit{trend filtering}~\citep{kim2009ell_1}, \textit{wavelet smoothing}~\citep{donoho1995adapting}, among others. 
   
When dealing with the generalized~\texttt{Lasso}, the inclusion of an $m \times d$ matrix $F$ in the total variation regularizer poses challenges for basic proximal gradient methods, such as the \textit{iterative shrinkage thresholding algorithm} (\texttt{ISTA}) and the \textit{fast iterative shrinkage thresholding algorithm} (\texttt{FISTA}).  A well-suited optimization method for solving the generalized \texttt{Lasso} is the~\textit{Alternating Direction Method of Multipliers} (\texttt{ADMM}), which was pioneered in the seventies of the last century by~\citet{glowinski1975approximation} and~\citet{gabay1976dual}. \texttt{ADMM} has gained popularity in the realm of machine learning due to its effectiveness in distributed convex optimization and its capability to handle large-scale optimization problems~\citep{boyd2011distributed}. However, the implementation of~\texttt{ADMM} faces the formidable task of computing the inverse of the matrix $A^{\top}A + sF^{\top}F$, where $s$ denotes the step size as indicated in~\citep{li2024understanding}. As its size increases, the task of inverting this matrix becomes more challenging. 

Before exploring saddle-point methods, let us first incorporate the generalized~\texttt{Lasso} into a general optimization problem represented as: 
\begin{equation}
\label{eqn: general-lasso}
\min_{x \in \mathbb{R}^{d_1}} \Phi(x) := f(x) + g(Fx),
\end{equation}
where both $f$ and $g$ are convex functions. By utilizing the conjugate transformation, also referred to as the Legendre-Fenchel transformation, we can reformulate the optimization problem~\eqref{eqn: general-lasso} in the following minimax form as
\begin{equation}
\label{eqn: game-rep}
\min_{y \in \mathbb{R}^{d_2}} \max_{x \in \mathbb{R}^{d_1}} \Phi(x,y) = \min_{x \in \mathbb{R}^{d_1}} \max_{y \in \mathbb{R}^{d_2}} \Phi(x,y) := f(x) + \big\langle Fx, y \big\rangle - g^{\star}(y),
\end{equation}
which was pioneered in~\citep{arrow1958studies}. The minimax form~\eqref{eqn: game-rep} allows us for the practical implementation of proximal operations. The proximal Arrow-Hurwicz algorithm, as introduced in~\citep{zhu2008efficient, esser2010general}, involves a descent iteration for the primal variable $x$ and an ascent iteration for the dual variable $y$ as follows: 
\begin{subequations}
\begin{empheq}[left=\empheqlbrace]{align}
  & x_{k+1} = \argmin_{x \in \mathbb{R}^{d_1}} \left\{ f(x) + \big\langle Fx, y_k \big\rangle + \frac{1}{2s} \|x - x_k \|^2 \right\},                         \label{eqn: ah-descent} \\
  & y_{k+1} = \argmax_{y \in \mathbb{R}^{d_2}} \left\{ - g^{\star}(y) + \big\langle Fx_{k+1}, y \big\rangle - \frac{1}{2s} \|y - y_k \|^2 \right\}.  \label{eqn: ah-ascent}
\end{empheq}
\end{subequations}
However, in certain scenarios,  the proximal Arrow-Hurwicz algorithm may fail to converge, as demonstrated by a counterexample in~\citep{he2014convergence}. To improve its practical effectiveness, incorporating a momentum step has been found to be advantageous. This enhanced algorithm, initially proposed in~\citep{chambolle2011first}, modifies the proximal Arrow-Hurwicz algorithm,~\eqref{eqn: ah-descent} and~\eqref{eqn: ah-ascent}, as
\begin{subequations}
\begin{empheq}[left=\empheqlbrace]{align}
  & x_{k+1} = \argmin_{x \in \mathbb{R}^{d_1}} \left\{ f(x) + \big\langle Fx, y_k \big\rangle + \frac{1}{2s} \|x - x_k \|^2 \right\},                                            \label{eqn: pdhg-descent} \\
  & \overline{x}_{k+1} = x_{k+1} + (x_{k+1} - x_k),                                                                                                                                                                        \label{eqn: pdhg-special} \\
  & y_{k+1} = \argmax_{y \in \mathbb{R}^{d_2}} \left\{ - g^{\star}(y) + \big \langle F\overline{x}_{k+1}, y \big\rangle - \frac{1}{2s} \|y - y_k \|^2 \right\},     \label{eqn: pdhg-ascent}
\end{empheq}
\end{subequations}
 which is also known as the~\textit{primal-dual hybrid gradient} (\texttt{PDHG}) algorithm. It is worth noting that the optimization problems addressable by the proximal Arrow-Hurwicz algorithm and~\texttt{PDHG} are equivalent to those solvable by~\texttt{ADMM}, as indicated in~\citep{li2024understanding}.  In general, the implicit solution in the iteration~\eqref{eqn: ah-ascent} or~\eqref{eqn: pdhg-ascent} is always directly obtainable when the convex function $g$ is either the $\ell_1$-norm or a quadratic function.  However, when the convex function $f$ is assumed to be $L$-smooth, finding the implicit solution in the iteration~\eqref{eqn: ah-descent} or~\eqref{eqn: pdhg-descent} requires the use of a gradient-based algorithm for approximation. Therefore, for the scope of this paper, we only focus on scenarios where f is either a quadratic function for the generalized \texttt{Lasso} or an indicator function for the~\textit{basis pursuit}~\citep{bruckstein2009sparse}. In such situations, the implicit solution is directly attainable.  Further theoretical insights on ergodic convergence can be found in~\citep{chambolle2016ergodic}.  The convergence behavior of~\texttt{PDHG}, from the contraction perspective, is explored in several studies,  including~\citep{he2012convergence, he2014convergence, he2022generalized}.

When considering a simple comparison between the two algorithms above,  the proximal Arrow-Hurwicz algorithm and \texttt{PDHG},  it naturally leads us to pose the following two questions:
\begin{tcolorbox}
\begin{itemize}
\item Why does~\texttt{PDHG} always converge on the minimax problem~\eqref{eqn: game-rep} whereas the proximal Arrow-Hurwicz algorithm may fail to converge in certain scenarios? 
\item What are the fundamental differences between these two algorithms?
\end{itemize}
\end{tcolorbox}
\noindent In this paper, we seek to explore these questions by leveraging techniques borrowed from ordinary differential equations (ODEs) and numerical analysis. Recent progress in these areas has notably advanced our understanding of the convergence behaviors of the vanilla gradient descent and shed light on the intriguing accelerated algorithms.  By utilizing Lyapunov analysis and phase-space representation, a comprehensive framework based on high-resolution ODEs has been well established in the studies~\citep{shi2022understanding, shi2019acceleration, chen2022gradient, chen2022revisiting, chen2023underdamped, li2023linear}.  Furthermore, this framework has also been extended to encompass proximal algorithms, including \texttt{ISTA} and \texttt{FISTA}s, in~\citep{li2022linear, li2022proximal, li2023linear}.  Most recently,  a new system of high-resolution ODEs has been derived to understand and analyze~\texttt{ADMM}, as presented in~\citep{li2024understanding}.  These advancements open up several new avenues for further research and development in optimization algorithms.

%
%
%
%
%
%
%

\subsection{$x$-correction and $y$-correction: a coupled small and essential perturbation}
\label{subsec: comparison-ah-pdhg}

When comparing the proximal Arrow-Hurwicz algorithm with~\texttt{PDHG} to study the iterative behavior of~\texttt{PDHG}, it is necessary to assume that both the objective functions, $f$ and $g$, are sufficiently smooth. This assumption allows us to eliminate the $\argmin$ and $\argmax$ operations to unfold the iterations. For the proximal Arrow-Hurwicz algorithm, the iterations,~\eqref{eqn: ah-descent} and~\eqref{eqn: ah-ascent}, can be expressed as
\begin{subequations}
\begin{empheq}[left=\empheqlbrace]{align}
  & \frac{x_{k+1} - x_{k}}{s} =  - F^{\top} y_k - \nabla f(x_{k+1}),                 \label{eqn: ah1-descent} \\
  & \frac{y_{k+1} - y_{k}}{s} = Fx_{k+1}  - \nabla g^{\star}(y_{k+1}).     \label{eqn: ah1-ascent}
\end{empheq}
\end{subequations}
If both functions, $f$ and $g$, are either constant or linear, the proximal Arrow-Hurwicz algorithm,~\eqref{eqn: ah1-descent} and~\eqref{eqn: ah1-ascent}, follows a forward-backward scheme and bears resemblance to the symplectic Euler scheme of a Hamilton system. Particularly, this resemblance is exact when the matrix $A$ is symmetric, aligning the proximal Arrow-Hurwicz algorithm,~\eqref{eqn: ah1-descent} and~\eqref{eqn: ah1-ascent}, directly with Hamiton systems. Further elaboration is provided in~\Cref{subsec: ah-non-convergence}, including the counterexample in~\citep{he2014convergence}.  Before moving on to explore~\texttt{PDHG}, we initiate by inserting $F^{\top} (y_{k+1} - y_{k})$ into both sides of the iteration~\eqref{eqn: ah1-descent}, transforming it into an implicit form as a preliminary step. By substituting~\eqref{eqn: pdhg-special} into~\eqref{eqn: pdhg-ascent}, we  expand the iteration~\eqref{eqn: pdhg-ascent} to express~\texttt{PDHG},~\eqref{eqn: pdhg-descent} ---~\eqref{eqn: pdhg-ascent}, as
\begin{subequations}
\begin{empheq}[left=\empheqlbrace]{align}
  & \frac{x_{k+1} - x_{k}}{s} - \underbrace{F^{\top} (y_{k+1} - y_{k})}_{y-\mathrm{correction}} = - F^{\top} y_{k+1} - \nabla f(x_{k+1}),        \label{eqn: pdhg1-descent} \\ 
  & \frac{y_{k+1} - y_{k}}{s} - \underbrace{F (x_{k+1} - x_{k})}_{x-\mathrm{correction}}             =   Fx_{k+1}   - \nabla g^{\star}(y_{k+1}),     \label{eqn: pdhg1-ascent}
\end{empheq}
\end{subequations}
where it can be observed upon the comparison of dimensions that $F^{\top} (y_{k+1} - y_{k})$ serves as a perturbation in the first iteration~\eqref{eqn: pdhg1-descent}, referred to as $y$-correction and $F (x_{k+1} - x_{k})$ serves as a perterbation in the second iteration~\eqref{eqn: pdhg1-ascent}, referred to as $x$-correction.  Importantly, the corrections, $x$-correction and $y$-correction, are not independent but rather interconnected, forming a coupled perturbation that crossed influences the iterations of~\texttt{PDHG},~\eqref{eqn: pdhg1-descent} and~\eqref{eqn: pdhg1-ascent}.  Finally, it is worth noting that the~\texttt{PDHG} algorithm,~\eqref{eqn: pdhg1-descent} and~\eqref{eqn: pdhg1-ascent}, also features the trait of an implicit scheme, characterized by a small but essential perturbation, which is akin to the scenario of~\texttt{ADMM} described in~\citep{li2024understanding}.

\subsection{Overview of contributions}
\label{subsec: overview}

In this paper, we employ techniques borrowed from ODEs and numerical analysis to understand and analyze the~\texttt{PDHG} algorithm and highlight our contributions as follows.

\begin{itemize}
\item[(1)] We utilize dimensional analysis, as proposed in~\citep{shi2022understanding, shi2023learning, shi2021hyperparameters, li2024understanding}, to establish a system of high-resolution ODEs for~\texttt{PDHG}, expressed as
\begin{subequations}
\begin{empheq}[left=\empheqlbrace]{align}
  & \dot{X} - s F^{\top} \dot{Y} =-  F^{T}Y - \nabla f(X)                  \label{eqn: high-descent} \\
  & \dot{Y} - s F\dot{X} = FX - \nabla g^{\star}(Y)                          \label{eqn: high-ascent}
\end{empheq}
\end{subequations}
It is important to note that the implicit (Euler) scheme of this system,~\eqref{eqn: high-descent} and~\eqref{eqn: high-ascent}, exactly corresponds to the~\texttt{PDHG} algorithm,~\eqref{eqn: pdhg-descent} ---~\eqref{eqn: pdhg-ascent}. Furthermore, the system of high-resolution ODEs,~\eqref{eqn: high-descent} and~\eqref{eqn: high-ascent}, effectively encapsulates the impact of coupled corrections, $x$-correction and $y$-correction. This marks a significant feature of~\texttt{PDHG}, distinguishing it from the proximal Arrow-Hurwicz algorithm~\eqref{eqn: ah-descent} and~\eqref{eqn: ah-ascent}.  The small but essential perturbation ensures that the variable pair $(X, Y)$ consistently converges, steering clear of the periodic behavior observed in the proximal Arrow-Hurwicz algorithm.

\item[(2)] Starting from the system of high-resolution ODEs,~\eqref{eqn: high-descent} and~\eqref{eqn: high-ascent}, we can devise a quadratic Lyapunov function and straightforwardly demonstrate that it decreases monotonically in the continuous scenario. This analysis seamlessly transitions from the continuous system of high-resolution ODEs,~\eqref{eqn: high-descent} and~\eqref{eqn: high-ascent}, to the discrete~\texttt{PDHG} algorithm,~\eqref{eqn: pdhg-descent} ---~\eqref{eqn: pdhg-ascent}. It is worth noting that the numerical error resulting from implicit discretization leads to the convergence of~\texttt{PDHG} with a rate of $O(1/N)$, echoing the phenomenon observed for the~\texttt{ADMM} algorithm as identified in~\citep{li2024understanding}. Furthermore, this analytical framework can be extended to the general form of~\texttt{PDHG} that involves two distinct parameters. Compared to earlier findings published in~\citep{chambolle2011first, chambolle2016ergodic, he2012convergence, he2014convergence, he2022generalized}, our proofs stand out for being principled, succinct, and straightforward.

\item[(3)] In addition, we further discover that if the objective function $f$ is strongly convex, the iterative sequence of~\texttt{PDHG} exhibits strong convergence in terms of average, with the convergence rate given by
\[
\|\overline{x}_N - x^{\star} \|^2 \leq O\left( \frac{1}{N+1} \right)
\]
where $x_N$ represents the average of the iterates, expressed as $1/N \sum_{k=1}^{N}x_k$. This assumption of strong convexity is applicable in practical scenarios, particularly for problems like the generalized \texttt{Lasso}.
\end{itemize}

%% file: 02_prelim.tex
\section{Preliminaries}
\label{sec: prelim}

In this section, we present a concise overview of basic definitions and classical theorems in convex optimization. These definitions and theorems are primarily drawn from classical literature, such as~\citep{rockafellar1970convex, rockafellar2009variational, nesterov1998introductory, boyd2004convex}, and will serve as references for proofs in the subsequent sections. Let us begin by defining a convex function, its conjugate function, and the concepts of subgradients and subdifferential. 
\begin{defn}
\label{defn: convex} 
A function $f: \mathbb{R}^{d} \mapsto \mathbb{R}$ is said to be convex if, for  any $x, y \in \mathbb{R}^d$ and any $\alpha \in [0,1]$, the following inequlaity holds:
\[
f\left( \alpha x + (1 - \alpha)y \right) \leq \alpha f(x) + (1 - \alpha)f(y).
\] 
The function $f^{\star}: \mathbb{R}^{d'} \mapsto \mathbb{R}$ is said to be the conjuate of the function $f$ if, for any $y \in \mathbb{R}^d$,  it satisfies the following definition as
\[
f^{\star}(y):= \sup_{x \in \mathbb{R}^{d}} \left( \big\langle y, x \big\rangle - g(x)\right). 
\]
\end{defn}

\begin{defn}
\label{defn: subgradient}
 A vector $v$ is said to be a subgradient of $f$ at $x$ if, for any $y \in \mathbb{R}^d$, the following inequlity holds:
\[
f(y) \geq f(x) + \langle v, y - x\rangle.
\] 
Additionally, the collection of all subgradients of $f$ at $x$ is referred to as the subdifferential of $f$ at $x$, and we denote it as $\partial f(x)$. 
\end{defn}

Let $\mathcal{C}(\mathbb{R}^d)$ be the class of continuous functions on $\mathbb{R}^d$. We can denote $\mathcal{F}^0(\mathbb{R}^d) \subset \mathcal{C}(\mathbb{R}^d)$ as the class of continuous and convex functions, and $\mathcal{F}^1(\mathbb{R}^d) \subset \mathcal{F}^0(\mathbb{R}^d)$ as its subclass, the class of differentiable and convex functions. According to~\Cref{defn: convex},  it is straightforward for us to conclude that for any function in $\mathcal{F}^1(\mathbb{R}^d)$, its conjugate belongs to $\mathcal{F}^1(\mathbb{R}^{d'})$. As stated in~\citep[ Proposition 8.21]{rockafellar1970convex}, since the function value is finite everywhere for any $f \in \mathcal{F}^0(\mathbb{R}^d)$,  subgradients always exist. This leads us to the following theorem. 

\begin{theorem}
\label{thm: subgradient-exist}
For any $f \in \mathcal{F}^{0}(\mathbb{R}^d)$, the subdifferential $\partial f(x)$ is nonempty for any $x \in \mathbb{R}^{d}$. 
\end{theorem}

Furthermore, we can establish a sufficient and necessary condition for a point to be a minimizer of any continuous and convex function. 

\begin{theorem}
\label{thm: minimizer}
A point $x^{\star}$ is a minimizer of $f \in \mathcal{F}^0(\mathbb{R}^d)$ if and only if $0 \in \partial f(x^{\star})$. 
\end{theorem}

Moving forward, we outline two basic properties of subdifferentials and subgradients, as collected from~\citep{rockafellar1970convex}. 

\begin{theorem}[Theorem 23.8 in \citep{rockafellar1970convex}]
\label{thm: sum-diff}
For any $f_1, f_2 \in\mathcal{F}^0(\mathbb{R}^d)$, the subdifferentials satisfy
\[
\partial(f_1 + f_2)(x) = \partial f_1(x) + \partial f_2(x), 
\]
for any $x \in \mathbb{R}^d$. 
\end{theorem}

\begin{theorem}[Theorem 25.1 in \citep{rockafellar1970convex}]
\label{thm: subgrad-unique}
For any $f \in \mathcal{F}^1(\mathbb{R}^d)$, the gradient $\nabla f(x)$ is the unique subgradient of $f$ at $x$, satisfying 
\[
f(y) \geq f(x) + \langle \nabla f(x), y - x\rangle,
\] 
for any $y \in \mathbb{R}^d$.
\end{theorem}

For smooth functions, we can provide precise characterizations as follows. 
\begin{defn}
\label{defn: L-smooth-strongly}
Let $f \in \mathcal{F}^1(\mathbb{R}^d)$. For any $x, y \in \mathbb{R}^d$, 
\begin{itemize}
\item[(1)] if the gradient of the function $f$ satisfies
\[
\|\nabla f(x) - \nabla f(y)\| \leq L \|x - y\|,
\]
we say that $f$ is $L$-smooth;
\item[(2)] if the function $f$ satisfies 
\[
f(y) \geq f(x) + \langle \nabla f(y), y - x \rangle + \frac{\mu}{2} \|y - x\|^2, 
\]
we say that $f$ is $\mu$-strongly convex. 
\end{itemize}
We denote $\mathcal{S}_{\mu,L}^1(\mathbb{R}^d) \subset \mathcal{F}^1(\mathbb{R}^d)$ as the class of $L$-smooth and $\mu$-strongly convex functions.
\end{defn}

Given that the implicit solution in the iteration~\eqref{eqn: ah-descent} or~\eqref{eqn: pdhg-descent} may not always be directly obtainable,  especially in scenarios where $f$ is $L$-smooth, it is necessary to specify the class of the objective function $f$ that allows for the direct obtainment of the implicit solution in these iterations. Let us denote this class as $\mathcal{R}(\mathbb{R}^d)$.  Finally, we introduce the definition of saddle points as denoted in~\citep[Definition 11.49]{rockafellar2009variational}. 
\begin{defn}
\label{defn: saddle}
 A vector pair $(x^{\star},y^{\star}) \in \mathbb{R}^{d_1} \times \mathbb{R}^{d_2}$ is said to be a saddle point of the convex-concave function $\Phi$ (convex with respect to the variable $x$ and concave with respect to the variable $y$) if, for any $x \in \mathbb{R}^{d_1}$ and $y \in \mathbb{R}^{d_2}$, the following inequality holds: 
\[
\Phi(x^{\star},y) \leq \Phi(x^{\star},y^{\star}) \leq  \Phi(x,y^{\star}).
\]
\end{defn}
With~\Cref{defn: saddle},  it ibecomes straightforward for us to derive the sufficient and necessary condition of the saddle point for the objective function $\Phi(x,y)$ taking the convex-concave form~\eqref{eqn: game-rep}, and we conclude this section with the following theorem.
\begin{theorem}
\label{thm: saddle-sn}
Let $f \in \mathcal{F}^1(\mathbb{R}^{d_1})$, $g \in \mathcal{F}^1(\mathbb{R}^{d_1})$ and $F \in \mathbb{R}^{d_2 \times d_1}$ be an $d_2 \times d_1$ matrix. Given the objective function $\Phi(x,y)$ takes the convex-concave form~\eqref{eqn: game-rep},  the point $(x^{\star},y^{\star})$ is a saddle if and only if the following inequalities holds
\[
\left\{\begin{aligned}
          & f(x) - f(x^{\star}) + \big\langle F(x - x^{\star}), y^{\star} \big\rangle \geq 0,                            \\ 
          & g^{\star}(y) - g^{\star}(y^{\star}) -  \big\langle  Fx^{\star}, y - y^{\star} \big\rangle \geq 0.         
         \end{aligned} \right.
\]          
\end{theorem}

%% file: 03_ode.tex
\section{Perspective from the system of high-resolution ODEs }
\label{sec: ode}

%
%

In this section, we first derive the system of high-resolution ODEs,~\eqref{eqn: high-descent} and~\eqref{eqn: high-ascent}, for the~\texttt{PDHG} algorithm through dimensional analysis. We then illustrate that the proximal Arrow-Hurwicz algorithm,~\eqref{eqn: ah1-descent} and~\eqref{eqn: ah1-ascent}, essentially acts as a numerical discretization for the system of low-resolution ODEs. Through the lens of Hamilton mechanics and its symplectic scheme that preserves the structure, we provide in-depth insights into the failure of the proximal Arrow-Hurwicz algorithm to converge in the counterexample given in~\citep{he2014convergence}, and we extend this analysis to encompass a range of counterexamples. Finally, we employ Lyapunov analysis to elucidate the convergence behaviors of the system of high-resolution ODEs.

\subsection{Derivation of the system of high-resolution ODEs}
\label{subsec: derivation-ode}

To derive the system of high-resolution ODEs, it is necessary to assume that both the objective functions and the solutions are sufficiently smooth. The parameter $s$ is set as the step size, which serves as a bridge between the~\texttt{PDHG} algorithm and the continuous system. Let $t_k = ks$, ($k=0,1,2,\ldots$) and take the ansatzs $x_k = X(t_k)$ and $y_k = Y(t_k)$. By performing a Taylor expansion in powers of $s$, we can obtain
\begin{subequations}
\begin{empheq}[left=\empheqlbrace]{align}
         x_{k} = X(t_{k}) & = x_{k+1} - s\dot{X}(t_{k+1}) + \frac{s^2}{2} \ddot{X}(t_{k+1}) + O\left(s^3\right),    \label{eqn: taylor-x}            \\
        y_{k} = Y(t_{k})  & = y_{k+1} - s\dot{Y}(t_{k+1})  + \frac{s^2}{2} \ddot{Y}(t_{k+1}) \;+ O\left(s^3\right). \label{eqn: taylor-y}   
\end{empheq}
\end{subequations}

For any $f \in \mathcal{F}^{1}(\mathbb{R}^d) \cap \mathcal{R}(\mathbb{R}^d)$ and $g \in \mathcal{F}^{1}(\mathbb{R}^d)$, we can expand~\texttt{PDHG} into the two iterations, denoted as ~\eqref{eqn: pdhg1-descent} and~\eqref{eqn: pdhg1-ascent}. Substituting~\eqref{eqn: taylor-x} into~\eqref{eqn: pdhg1-descent} and~\eqref{eqn: taylor-y} into~\eqref{eqn: pdhg1-ascent}, we have
\begin{subequations}
\begin{empheq}[left=\empheqlbrace]{align}
   &- \frac{s}{2} \ddot{X}(t_{k+1}) +  \dot{X}(t_{k+1}) - s F^{\top} \dot{Y}(t_{k+1}) +O(s^2)  = - F^{T}Y(t_{k+1}) - \nabla f(X(t_{k+1})),                      \label{eqn: high-descent-derive} \\
   &- \frac{s}{2} \ddot{Y}(t_{k+1}) +  \dot{Y}(t_{k+1}) - s F\dot{X}(t_{k+1}) + O(s^2) =  FX(t_{k+1}) - \nabla g^{\star}(Y(t_{k+1})).                              \label{eqn: high-ascent-derive}
\end{empheq}
\end{subequations}
By focusing the $O(s)$ terms and disregarding the $O(s^2)$ terms, we can derive a system of ODEs as
\[
\left\{ \begin{aligned}
            &- \frac{s}{2} \ddot{X} +  \dot{X} - s F^{\top} \dot{Y}  = - F^{T}Y - \nabla f(X),        \\
            &- \frac{s}{2} \ddot{Y} +  \dot{Y} - s F\dot{X} \; = FX - \nabla g^{\star}(Y),                            
          \end{aligned} \right. 
\]
which corresponds to an overdamped system since the step size $s$ is small. In an overdamped system, any inertial effects are damped out quickly by the effects of viscous, frictional, or other damping forces, allowing their effects to be disregarded, as noted in~\citep{adams2013large}. This concept is further exemplified in~\citep[Section 1.2]{shi2021hyperparameters}. Therefore, the system of high-resolution ODEs,~\eqref{eqn: high-descent} and~\eqref{eqn: high-ascent}, serves as a reasonable approximation for~\texttt{PDHG},~\eqref{eqn: pdhg1-descent} and~\eqref{eqn: pdhg1-ascent}.  In turn, the~\texttt{PDHG} algorithm,~\eqref{eqn: pdhg1-descent} and~\eqref{eqn: pdhg1-ascent}, exactly corresponds to  the implicit discretization of the system of high-resolution ODEs,~\eqref{eqn: high-descent} and~\eqref{eqn: high-ascent}.

\subsection{The system of low-resolution ODEs: Symplectic-Euler scheme and proximal Arrow-Hurwicz algorithm}
\label{subsec: ah-non-convergence}

By focusing the $O(1)$ terms and disregarding the $O(s)$ terms in equations~\eqref{eqn: high-descent-derive} and~\eqref{eqn: high-ascent-derive}, we can derive a simplified system as
\begin{subequations}
\begin{empheq}[left=\empheqlbrace]{align}
  & \dot{X} = - F^{T}Y - \nabla f(X),                               \label{eqn: low-descent} \\
  & \dot{Y}  =  FX - \nabla g^{\star}(Y),                          \label{eqn: low-ascent}
\end{empheq}
\end{subequations}
which is referred to as the system of low-resolution ODEs for~\texttt{PDHG}, in contrast to the system of high-resolution ODEs given by equations~\eqref{eqn: high-descent} and~\eqref{eqn: high-ascent}.  Furthermore, a mixed explicit-implicit Euler discretization of the system of low-resolution ODEs,~\eqref{eqn: low-descent} and~\eqref{eqn: low-ascent}, exactly corresponds with the proximal Arrow-Hurwicz algorithm,~\eqref{eqn: ah1-descent} and~\eqref{eqn: ah1-ascent}.

\begin{figure}[htb!]
\centering
\includegraphics[scale=0.28]{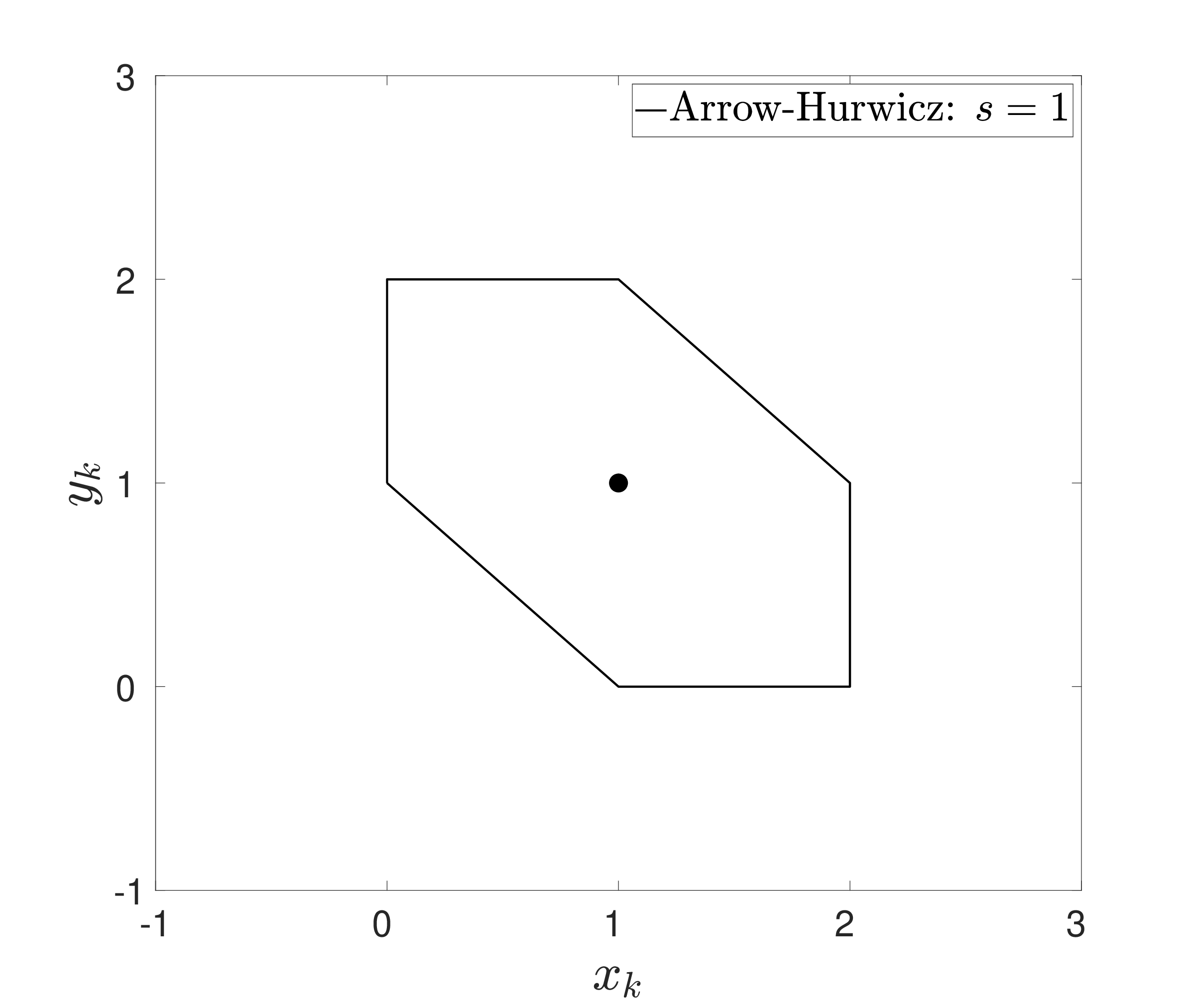}
\caption{The counterexample demonstrated in~\citep{he2014convergence}: given the objective function~\eqref{eqn: counter-pot}, the trajectory generated the proximal Arrow-Hurwicz algorithm,~\eqref{eqn: ah1-descent} and~\eqref{eqn: ah1-ascent}, starting from the point $(0,1)$ under the object, fails to converge to the saddle point $(1,1)$. }
\label{fig: s1}
\end{figure}

As stated in~\citep{he2014convergence},  the objective function is smooth and expressed as 
\begin{equation}
\label{eqn: counter-pot}
\Phi(x,y) = x - xy + y,
\end{equation} 
where we can observe that $(x^{\star}, y^{\star}) = (1,1)$ is the unique saddle point.  Substituting the objective function~\eqref{eqn: counter-pot} into the proximal Arrow-Hurwicz algorithm,~\eqref{eqn: ah1-descent} and~\eqref{eqn: ah1-ascent}, yields
\begin{subequations}
\begin{empheq}[left=\empheqlbrace]{align}
          & \frac{x_{k+1} - x_{k}}{s} =  y_{k} - 1,          \label{eqn: sym-eu-x}         \\
          & \frac{y_{k+1} - y_{k}}{s} = - x_{k+1} +1.    \label{eqn: sym-eu-y}
\end{empheq}
\end{subequations}
When the step size is set to $s = 1$, equations~\eqref{eqn: sym-eu-x} and~\eqref{eqn: sym-eu-y} correspond to the counterexample demonstrated in~\citep{he2014convergence}. This scenario is illustrated in~\Cref{fig: s1} where the iterative trajectory of the proximal Arrow-Hurwicz algorithm,~\eqref{eqn: ah1-descent} and~\eqref{eqn: ah1-ascent}, forms a closed loop consisting of six points, failing to converge to the unique saddle point $(x^{\star}, y^{\star}) = (1,1)$ . 
Furthermore, equations~\eqref{eqn: sym-eu-x} and~\eqref{eqn: sym-eu-y} can be interpreted as the symplectic Euler scheme of the following Hamilton system: 
\begin{subequations}
\begin{empheq}[left=\empheqlbrace]{align}
          & \dot{X} =  Y - 1,              \label{eqn: ham-x} \\
          & \dot{Y} =  - X +1.           \label{eqn: ham-y}   
\end{empheq}
\end{subequations}
It can be observed that for the given objective function~\eqref{eqn: counter-pot}, the system of low-resolution ODEs,~\eqref{eqn: low-descent} and~\eqref{eqn: low-ascent}, exactly matches the Hamilton system,~\eqref{eqn: ham-x} and~\eqref{eqn: ham-y}.

\begin{figure}[htb!]
\centering
\includegraphics[scale=0.28]{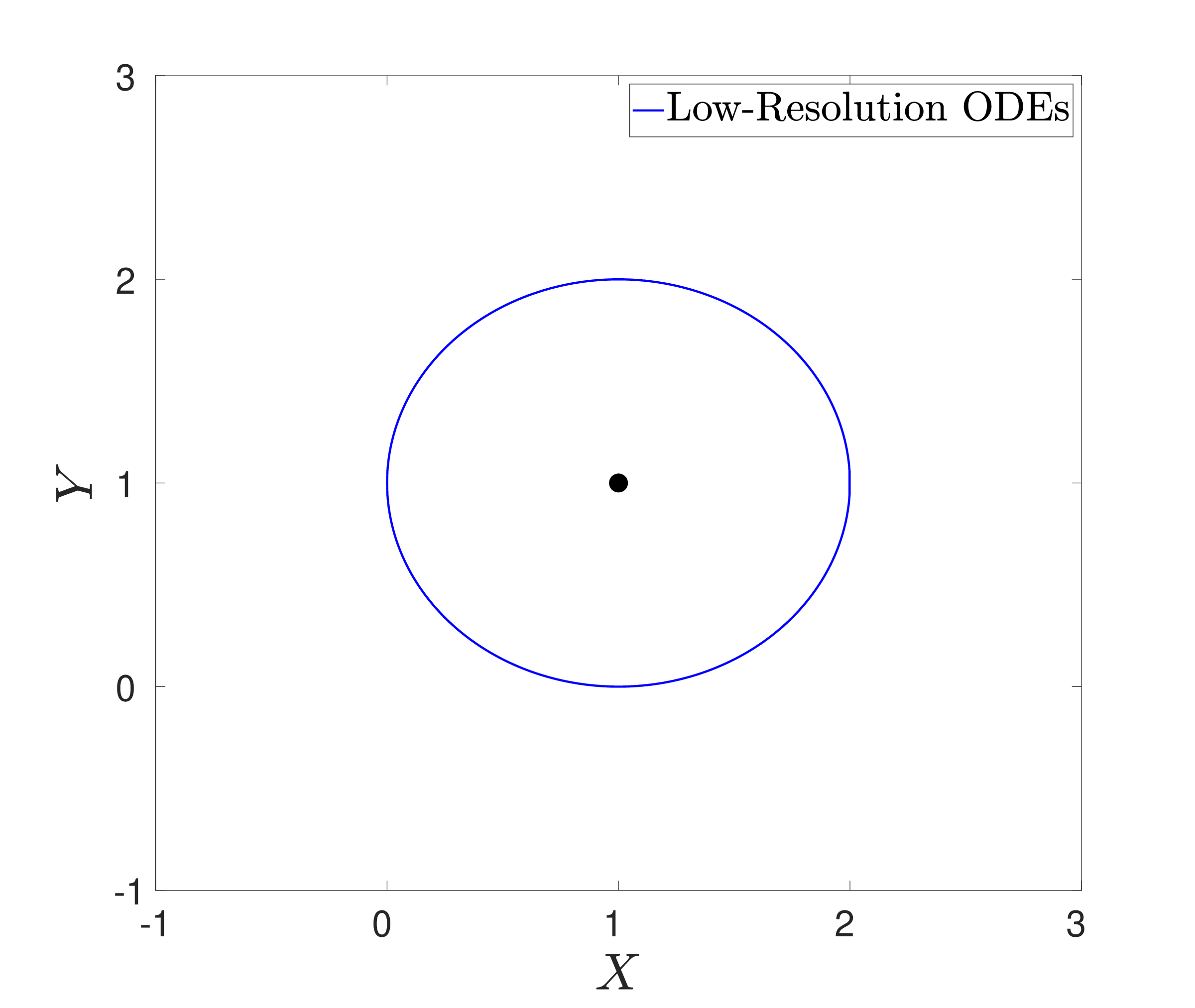}
\caption{Given the objective function~\eqref{eqn: counter-pot}, the trajectory of the system of low-resolution ODEs,~\eqref{eqn: low-descent} and~\eqref{eqn: low-ascent}, starting from the point $(0,1)$ .}
\label{fig: ode}
\end{figure}

\begin{figure}[htb!]
\centering
\begin{subfigure}[t]{0.45\linewidth}
\centering
\includegraphics[scale=0.18]{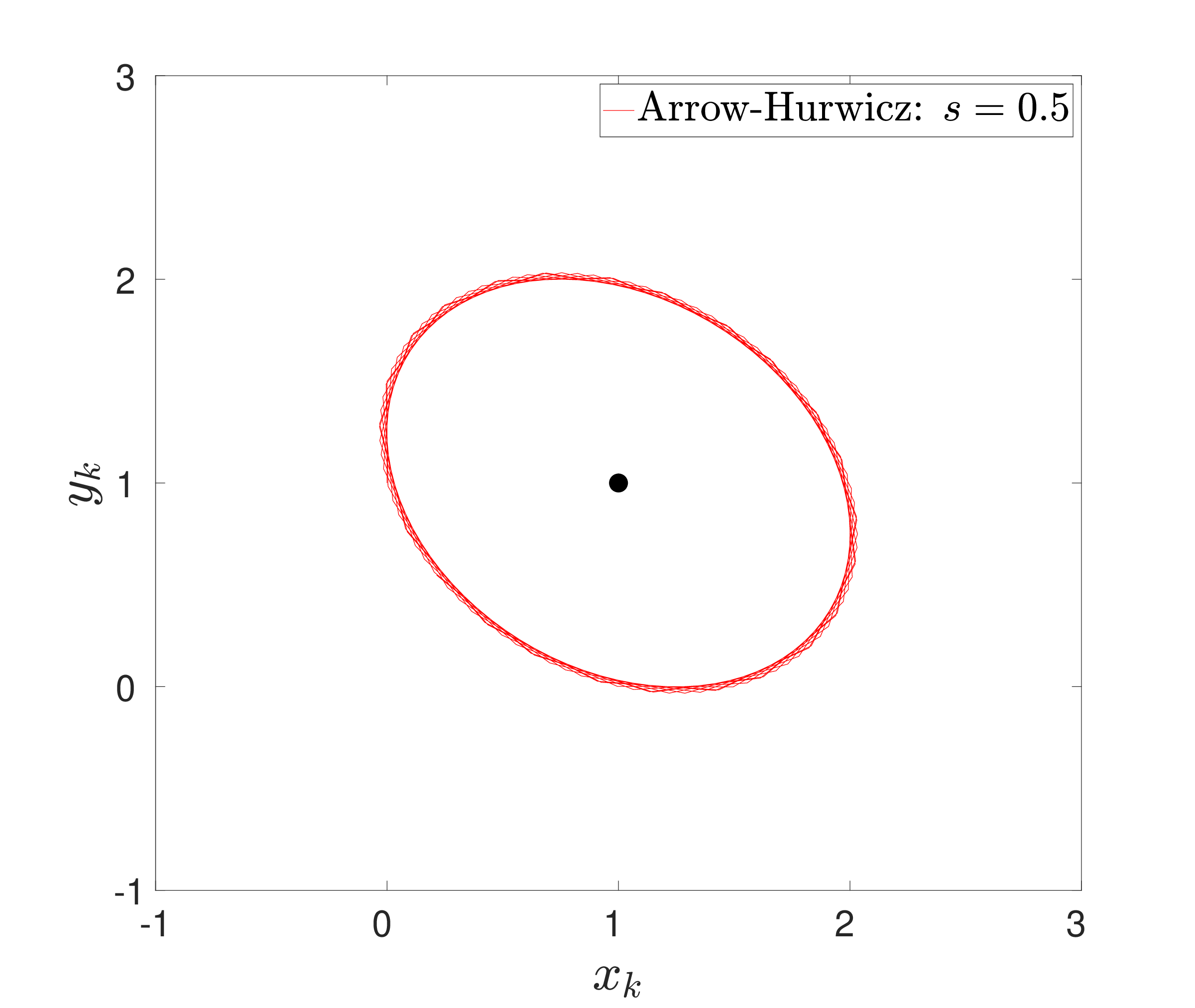}
\caption{$s=0.5$}
\label{subfig: s05}
\end{subfigure}
\begin{subfigure}[t]{0.45\linewidth}
\centering
\includegraphics[scale=0.18]{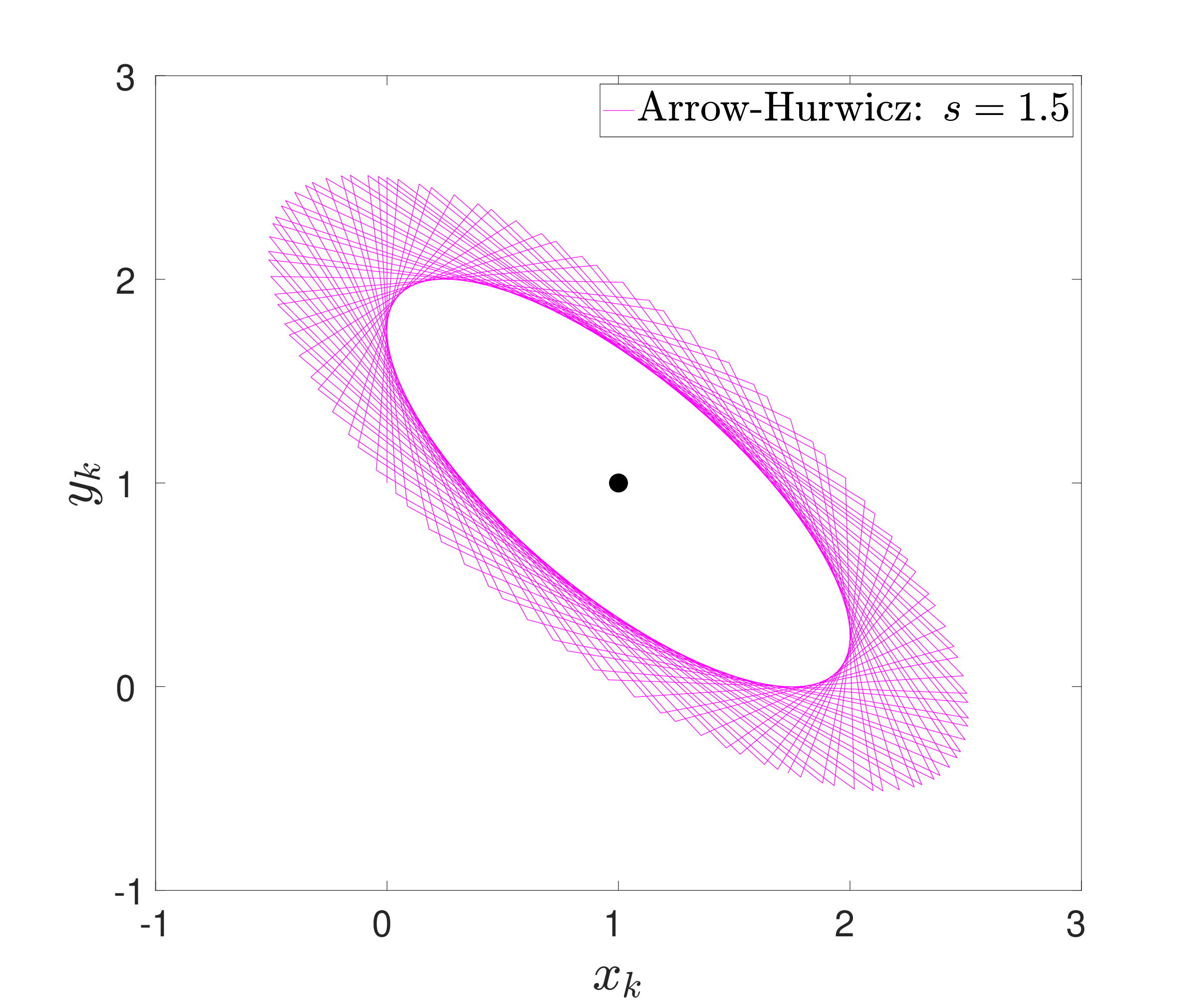}
\caption{$s=1.5$}
\label{subfig: s15}
\end{subfigure}
\caption{Given the objective function~\eqref{eqn: counter-pot}, the trajectories generated the proximal Arrow-Hurwicz algorithm,~\eqref{eqn: ah1-descent} and~\eqref{eqn: ah1-ascent}, starting from $(0,1)$ with different step sizes.}
\label{fig: s0515}
\end{figure}

By analyzing the Hamilton dynamics outlined in equations~\eqref{eqn: ham-x} and~\eqref{eqn: ham-y}, we can gain more insight into why the proximal Arrow-Hurwicz algorithm fails to converge from the perspective of Hamilton mechanics. As demonstrated in~\Cref{fig: ode}, the trajectory of the Hamilton system, as governed by equations~\eqref{eqn: ham-x} and~\eqref{eqn: ham-y}, originates from the point $(0,1)$ and traces a perfect circle centered at $(1,1)$ with a radius $r=1$. The dynamic behavior is encapsulated by the Hamilton function 
\[
H(x,y) = \frac{1}{2}(x^2 + y^2) - x - y = - \frac12.
\]
Essentially, the Hamilton system keeps the value of the Hamilton function constant. At the starting point $(0,1)$, the Hamilton function holds a value $H = - 1/2$  while at the unique saddle $(1,1)$, the value shifts to $H = -1$. Consequently, the Hamilton system governed by equations~\eqref{eqn: ham-x} and~\eqref{eqn: ham-y} fails to converge towards the saddle point. In other words, given the objective function~\eqref{eqn: counter-pot}, the system of low-resolution ODEs,~\eqref{eqn: low-descent} and~\eqref{eqn: low-ascent}, does not approach the saddle point. Even though numerical errors might lead to deviations in the exact preservation of the value of the Hamilton function,  the periodic structure can still be sustained in a topological sense. The symplectic Euler scheme, as elaborated in~\citep{feng1984difference, hairer2006geometric}, employs both a forward iteration (e.g.~\eqref{eqn: sym-eu-x}) and a backward iteration (e.g.~\eqref{eqn: sym-eu-y}) to preserve this structure effectively. Through an analysis of the discriminant of a quadratic equation, it is determined that the periodic structure can be maintained under the symplectic Euler scheme,~\eqref{eqn: sym-eu-x} and~\eqref{eqn: ham-y}, provided the step size falls within the range $0 < s < 2$. Building on the explanation above, it can be inferred that the proximal Arrow-Hurwicz algorithm, as outlined in~\eqref{eqn: ah1-descent} and~\eqref{eqn: ah1-ascent}, always maintain the periodic structure when applied to the objective function~\eqref{eqn: counter-pot}. As a result, it cannot converge towards the saddle point $(1, 1)$ for any step size within the range $0 < s < 2$. Furthermore, it is noteworthy that  the failure of the proximal Arrow-Hurwicz algorithm to converge is not limited to a single counterexample, as pointed out in~\citep{he2014convergence}, but extends to a series of counterexamples. These are depicted through two examples in~\Cref{fig: s0515}, highlighting the broader implications and the diverse scenarios where the proximal Arrow-Hurwicz algorithm preserves the periodic structure and falls short of achieving convergence.


\subsection{Convergence of the system of high-resolution ODEs}
\label{subsec: convergence-ode}

Let $(X,Y)$ be the solution to the system of high-resolution ODEs,~\eqref{eqn: high-descent} and~\eqref{eqn: high-ascent}, and $(x,y)$ be any point in $\mathbb{R}^{d_1} \times \mathbb{R}^{d_2}$.  To analyze the convergence of the system of high-resolution ODEs, we construct the following Lyapunov function as
\begin{equation}
\label{eqn: lyapunov-ode}
\mathcal{E}(t) = \frac{1}{2s} \|X - x\|^2 + \frac{1}{2s} \|Y - y\|^2 - \big\langle F(X - x), Y - y \big\rangle,
\end{equation}
where the Lyapunov function $\mathcal{E}(t) \geq 0$ always holds for any step size satisfying $ 0 < s \|F\| \leq  1$.\footnote{ In this context, the matrix norm $\|\cdot\|$ refers to the spectral norm, which is induced by the $\ell_2$ norm for vectors, given as \[ \|F\| = \sqrt{\lambda_{\max}(F^{\top}F)}.\]} 
\begin{lemma}
\label{lem: high-ode}
Let $f \in \mathcal{F}^{1}(\mathbb{R}^{d_1}) \cap \mathcal{R}(\mathbb{R}^{d_1}) $ and $g \in \mathcal{F}^{1}(\mathbb{R}^{d_1})$. For any step size $s \in \big(0, \|F\|^{-1}\big]$,  the Lyapunov function~\eqref{eqn: lyapunov-ode} satisfies the following inequality as
\begin{equation}
\label{eqn: high-ode-inq}
\frac{d\mathcal{E}}{dt} \leq  - \frac{1}{s} \big\langle  F(X - x), y  \big\rangle  + \frac{1}{s} \big\langle Fx , Y - y  \big\rangle + \frac{1}{s}\left( f(x) - f(X) + g^{\star}(y) - g^{\star}(Y) \right).
\end{equation}
\end{lemma}

\begin{proof}[Proof of~\Cref{lem: high-ode}]
With the system of high-resolution ODEs,~\eqref{eqn: high-descent} and~\eqref{eqn: high-ascent}, we can calculate  the time derivative of the Lyapunov function~\eqref{eqn: lyapunov-ode} as
\begin{align}
\frac{d\mathcal{E}}{dt} & = \frac{1}{s} \big\langle \dot{X}, X - x \big\rangle - \big\langle F^{\top}\dot{Y} , X - x\big \rangle + \frac{1}{s}\big \langle \dot{Y}, Y - y \big\rangle - \big\langle F\dot{X}, Y - y \big\rangle, \nonumber \\
                                     & =  - \frac{1}{s} \big\langle X - x, F^{\top}Y  \big\rangle  + \frac{1}{s} \big\langle FX , Y - y  \big\rangle - \frac{1}{s} \big\langle \nabla f(X), X - x  \big\rangle - \frac{1}{s} \big\langle \nabla g^{\star}(Y), Y - y  \big\rangle \nonumber \\
                                     & =   - \frac{1}{s} \big\langle F(X - x), y \big\rangle + \frac{1}{s} \big\langle Fx , Y - y  \big\rangle + \frac{1}{s} \big\langle \nabla f(X), x - X  \big\rangle + \frac{1}{s} \big\langle \nabla g^{\star}(Y),  y - Y \big\rangle, \label{eqn: lem-high-1}
\end{align}
where the last equality follows $\big\langle  X - x, F^{\top}(Y-y)  \big\rangle = \big\langle F(X - x), Y-y \big\rangle $. Since $f \in \mathcal{F}^{1}(\mathbb{R}^{d_1})$ and $g \in \mathcal{F}^1(\mathbb{R}^{d_1})$, we can derive the following convex inequalities from~\Cref{defn: convex} and~\Cref{thm: subgrad-unique} as
\begin{subequations}
\begin{empheq}[left=\empheqlbrace]{align}
          & f(x) - f(X)  \geq  \big\langle \nabla f(X), x - X  \big\rangle,                                             \label{eqn: lem-high-2f} \\
          &g^{\star}(y) - g^{\star}(Y) \geq  \big\langle \nabla g^{\star}(Y),  y - Y \big\rangle.           \label{eqn: lem-high-2g}   
\end{empheq}
\end{subequations}
By substituting~\eqref{eqn: lem-high-2f} and~\eqref{eqn: lem-high-2g} into~\eqref{eqn: lem-high-1}, we complete the proof. 
\end{proof}

In the given context,  the time average of a variable, denoted as $X$, within the time interval $[0,t]$, is defined as
\[
\overline{X} = \frac{1}{t} \int_{0}^{t} X(s)ds.
\]
We introduce the set of points that satisfies the variational inequality, denoted as $\mathrm{VI}$, as follows:
\begin{equation}
\label{eqn: vi}
\mathrm{VI} = \left\{ (x^{\diamond}, y^{\diamond}) \big| f(x^{\diamond}) - f(x) + g^{\star}(y^{\diamond}) - g^{\star}(y) + \big\langle F(x^{\diamond} - x), y \big\rangle  -  \big\langle Fx, y^{\diamond} - y \big\rangle \leq 0 \right\}.
\end{equation}
Utilizing~\Cref{lem: high-ode}, we can deduce the convergence rate of the time average to the set $\mathrm{VI}$ in a weak sense, which is rigorously stated in the theorem below.

\begin{theorem}
\label{thm: high-ode-weak}
Under the assumptions of~\Cref{lem: high-ode}, the time average $(\overline{X}, \overline{Y})$ converges weakly to the set $\mathrm{VI}$, denoted as~\eqref{eqn: vi}, at the  rate characterized by
\begin{multline}
\label{eqn: high-ode-ergodic}
f(\overline{X}) - f(x) + g^{\star}(Y) - g^{\star}(y) + \big\langle F(\overline{X} - x), y \big\rangle  -  \big\langle Fx, \overline{Y} - y \big\rangle \\ \leq \frac{\|x_0 - x\|^2 + \|y_0 - y\|^2 - 2s\big\langle F(x_0 - x), y_0 - y \big\rangle}{2t}
\end{multline} 
for any initial $(x_0, y_0) \in \mathbb{R}^{d_1} \times \mathbb{R}^{d_2}$. 
\end{theorem}

When the point $(x, y)$ is designated as the saddle point $(x^{\star},y^{\star})$, the Lyapunov function~\eqref{eqn: lyapunov-ode} becomes 
\begin{equation}
\label{eqn: lyapunov-ode-saddle}
\mathcal{E}(t) = \frac{1}{2s} \|X - x^{\star}\|^2 + \frac{1}{2s} \|Y - y^{\star}\|^2 + \big\langle F(X - x^{\star}), Y - y^{\star} \big\rangle.
\end{equation}
Substituting the saddle point $(x^{\star}, y^{\star})$ into the derivative inequality~\eqref{eqn: lem-high-1} yields
\begin{equation}
\label{eqn: high-ode-inq-saddle}
\frac{d\mathcal{E}}{dt} \leq  - \frac{1}{s} \big\langle F^{\top}y^{\star}, X - x^{\star}  \big\rangle  + \frac{1}{s} \big\langle Fx^{\star} , Y - y^{\star}  \big\rangle + \frac{1}{s}\left( f(x^{\star}) - f(X) + g^{\star}(y^{\star}) - g^{\star}(Y) \right).
\end{equation}
With the help of~\Cref{thm: saddle-sn}, it can be established that the Lyapunov function decreases monotonically. The rigorous representation is stated as follows.
\begin{theorem}
\label{thm: high-ode-ergodic}
Under the assumptions of~\Cref{lem: high-ode}, the Lyapunov function given by~\eqref{eqn: lyapunov-ode-saddle} decreases monotonically.
\end{theorem}

In addition, considering the generalized~\texttt{Lasso}, the objective function $f$ may be $\mu$-strongly convex. Assuming that the objective function $f$ does indeed possess $\mu$-strongly convexity, we can deduce that the time average of $X$ over the time interval $[0,t]$ converges strongly with a rate of $O(1/t)$.  The formal statement is given as follows.
\begin{theorem}
\label{thm: high-ode-strong}
Under the assumptions of~\Cref{lem: high-ode},  if we further assume the objective function satisfies $f \in \mathcal{S}_{\mu,L}^{1}(\mathbb{R}^d)$, the time average $\overline{X}$ converges at the rate characterized by
\begin{equation}
\label{eqn: high-ode-strong}
\big\| \overline{X} - x^{\star} \big\|^2 \leq \frac{\|x_0 - x^{\star}\|^2 + \|y_0 - y^{\star}\|^2 + 2s\big\langle F(X - x^{\star}), Y - y^{\star} \big\rangle}{\mu t}
\end{equation}
\end{theorem}
\begin{proof}[Proof of~\Cref{thm: high-ode-strong}]
Given the objective function is assumed to be $f \in \mathcal{S}_{\mu,L}^{1}(\mathbb{R}^d)$, we can refine the convex inequality~\eqref{eqn: lem-high-2f} as
\begin{equation}
\label{eqn: strong-inq}
 f(x) - f(X)  \geq  \big\langle \nabla f(X), x - X  \big\rangle + \frac{\mu}{2} \|x - X\|^2.
\end{equation}
By utilizing the refined inequality~\eqref{eqn: strong-inq} and following the same process that leads to~\Cref{thm: high-ode-ergodic}, it can be determined that the time derivative of the Lyapunov function~\eqref{eqn: lyapunov-ode-saddle} satisfies the following inequality:
\[
\frac{d \mathcal{E}}{dt} \leq - \frac{\mu}{2s} \|X - x^{\star}\|^2.
\]
Hence, by taking the time average of $X$ and leveraging the convexity of the $\ell_2$-norm square, we complete the proof.
\end{proof}


%
%

%% file: 04_algorithm.tex
\section{Convergence rates of PDHG}
\label{sec: pdhg}

In~\Cref{subsec: convergence-ode}, we have delved into the convergence behavior of the system of high-resolution ODEs,~\eqref{eqn: high-descent} and~\eqref{eqn: high-ascent}. Now, we extend our exploration straightforwardly from the continuous perspective to the discrete setting, shifting our focus on the~\texttt{PDHG} algorithm,~\eqref{eqn: pdhg-descent} ---~\eqref{eqn: pdhg-ascent}. It is important to highlight that the~\texttt{PDHG} algorithm,~\eqref{eqn: pdhg-descent} ---~\eqref{eqn: pdhg-ascent}, exhibits a phenomenon previously observed with~\texttt{ADMM}, as identified in~\citep[Section 1.2]{li2024understanding}, where numerical errors arising from implicit discretization have a significant impact on the convergence rate.

Consider any point $(x,y) \in \mathbb{R}^{d_1} \times \mathbb{R}^{d_2}$ and the iterative sequence $\{(x_k, y_k)\}_{k=0}^{\infty}$ generated by~\texttt{PDHG},~\eqref{eqn: pdhg-descent} ---~\eqref{eqn: pdhg-ascent}.  We extend the continuous Lyapunov function~\eqref{eqn: lyapunov-ode} to its discrete form counterpart as 
\begin{equation}
\label{eqn: pdhg-lyapunov}
\mathcal{E}(k) = \frac{1}{2s} \|x_k - x\|^2 + \frac{1}{2s} \|y_k - y\|^2 - \big\langle F(x_k - x), y_k - y \big\rangle,
\end{equation}
where the Lyapunov function $\mathcal{E}(k) \geq 0$ always holds for any step size satisfying $ 0 < s \|F\| \leq  1$. 
\begin{lemma}
\label{lem: pdhg}
Let $f \in \mathcal{F}^{0}(\mathbb{R}^{d_1}) \cap \mathcal{R}(\mathbb{R}^{d_1}) $ and $g \in \mathcal{F}^{0}(\mathbb{R}^{d_1})$. For any step size $s \in \big(0, \|F\|^{-1}\big]$, the discrete Lyapunov function~\eqref{eqn: pdhg-lyapunov} satisfies the following inequality as
\begin{align}
\mathcal{E}(k+1) - &\mathcal{E}(k) \nonumber \\
                   \leq & -  \big\langle F(x_{k+1} - x), y\big\rangle  + \big\langle Fx , y_{k+1} - y  \big\rangle +  f(x) - f(x_{k+1}) + g^{\star}(y) - g^{\star}(y_{k+1})  \nonumber \\
                           & - \underbrace{ \left(   \frac{1}{2s} \|x_{k+1} - x_{k}\|^2 +  \frac{1}{2s}  \|y_{k+1} - y_{k}\|^2 -  \big\langle F(x_{k+1} - x_{k}), y_{k+1} - y_{k} \big\rangle \right)}_{\textbf{NE}}, \label{eqn: pdhg-inq}
\end{align}
where $\mathbf{NE}$ represents the numerical error resulting from the implicit scheme.
\end{lemma}
\begin{proof}[Proof of~\Cref{lem: pdhg}]
We begin by examining the discrete Lyapunov function~\eqref{eqn: pdhg-lyapunov}. Our goal is to compute the difference in its value between successive iterations, i.e., $\mathcal{E}(k+1) - \mathcal{E}(k)$.  This difference can be expressed as: 
\begin{align}
    \mathcal{E}(k+1) - &\mathcal{E}(k) \nonumber  \\
= & \underbrace{\left\langle \frac{x_{k+1} - x_{k}}{s} - F^{\top}(y_{k+1} - y_{k}), x_{k+1} - x \right\rangle + \left\langle \frac{y_{k+1} - y_{k}}{s} - F(x_{k+1} - x_{k}), y_{k+1} - y \right\rangle}_{\mathbf{I}} \nonumber \\ 
                                                          & -  \underbrace{ \left( \frac{1}{2s} \|x_{k+1} - x_k\|^2 + \frac{1}{2s} \|y_{k+1} - y_{k}\|^2 - \big\langle F(x_{k+1} - x_{k}), y_{k+1} - y_{k} \big\rangle \right) }_{\mathbf{NE}}, \label{eqn: pdhg-lem1}
\end{align}
where $\mathbf{I}$ corresponds to the term derived from the continuous perspective, and $\mathbf{NE}$ represents the numerical error resulting from the implicit discretization.  Based on the identity $\big \langle x_{k+1} - x, F^{\top}(y_{k+1} - y)\big\rangle = \big \langle F(x_{k+1} - x), y_{k+1} - y\big\rangle$, we can reformulate $\mathbf{I}$ as
\begin{align}
\mathbf{I} = & \left\langle \frac{x_{k+1} - x_{k}}{s} - F^{\top}(y_{k+1} - y_{k}) + F^{\top}y_{k+1}, x_{k+1} - x \right\rangle   \nonumber \\
                    & + \left\langle \frac{y_{k+1} - y_{k}}{s} - F(x_{k+1} - x_{k}) - Fx_{k+1}, y_{k+1} - y \right\rangle  \nonumber \\ 
                    &- \big \langle F(x_{k+1} - x), y \big\rangle + \big\langle Fx, y_{k+1} - y \big\rangle. \label{eqn: pdhg-lem2}
\end{align}
Given any $f \in \mathcal{F}^{0}(\mathbb{R}^{d_1}) \cap \mathcal{R}(\mathbb{R}^{d_1}) $ and any $g \in \mathcal{F}^{0}(\mathbb{R}^{d_1})$,~\Cref{thm: minimizer} and~\Cref{thm: sum-diff} allow us to bypass the $\argmin$ and $\argmax$ operations and further explore~\texttt{PDHG} in terms of the first decent iteration~\eqref{eqn: pdhg-descent} and the third ascent iteration~\eqref{eqn: pdhg-ascent}, expressed as:
\begin{subequations}
\begin{empheq}[left=\empheqlbrace]{align}
          &  \frac{x_{k+1} - x_{k}}{s} + F^{\top}y_{k} + \partial f(x_{k+1}) \ni 0,                                \label{eqn: subgrad-pdhg-descent} \\
          &  \frac{y_{k+1} - y_{k} }{s} - F\overline{x}_{k+1} + \partial g^{\star}(y_{k+1}) \ni 0.          \label{eqn: subgrad-pdhg-ascent}   
\end{empheq}
\end{subequations}
By isolating the term $F^{\top}(y_{k+1} - y_{k})$ in the $x$-update~\eqref{eqn: subgrad-pdhg-descent} and incorporating the momentum step~\eqref{eqn: pdhg-special} into the $y$-update~\eqref{eqn: subgrad-pdhg-ascent}, we have
\begin{subequations}
\begin{empheq}[left=\empheqlbrace]{align}
          & \frac{x_{k+1} - x_{k}}{s} - F^{\top}(y_{k+1} - y_{k}) + F^{\top}y_{k+1}  + \partial f(x_{k+1})   \ni 0,                                \label{eqn: subgrad1-pdhg-descent} \\
          & \frac{y_{k+1} - y_{k}}{s} - F(x_{k+1} - x_{k}) - Fx_{k+1} + \partial g^{\star}(y_{k+1})    \ni 0.                                         \label{eqn: subgrad1-pdhg-ascent}   
\end{empheq}
\end{subequations}
Leveraging these updates,=, we can derive two convex inequalities based on~\Cref{defn: subgradient} and~\Cref{thm: sum-diff} as
\begin{subequations}
\begin{empheq}[left=\empheqlbrace]{align}
        &   f(x) - f(x_{k+1}) \geq \left\langle \frac{x_{k+1} - x_{k}}{s} - F^{T}(y_{k+1} - y_{k}) + F^{\top}y_{k+1}, x_{k+1} - x\right\rangle,         \label{eqn: subgrad2-pdhg-descent} \\
        &  g^{\star}(y) - g^{\star}(y_{k+1}) \geq  \left\langle \frac{y_{k+1} - y_{k}}{s} - F(x_{k+1} - x_{k}) - Fx_{k+1} , y_{k+1} - y \right\rangle.\label{eqn: subgrad2-pdhg-ascent}   
\end{empheq}
\end{subequations}
Substituting these two convex inequalities~\eqref{eqn: subgrad2-pdhg-descent} and~\eqref{eqn: subgrad2-pdhg-ascent} into our earlier formulation~\eqref{eqn: pdhg-lem2} allows us to establish an upper bound for the iterative difference:
\begin{equation}
\label{eqn: pdhg-lem3}
\mathbf{I} \leq - \big \langle F(x_{k+1} - x), y \big\rangle + \big\langle Fx_{k+1}, y_{k+1} - y \big\rangle + f(x) - f(x_{k+1}) + g^{\star}(y) - g^{\star}(y_{k+1}).  
\end{equation}
Finally, by inserting this inequality~\eqref{eqn: pdhg-lem3} into the iterative difference~\eqref{eqn: pdhg-inq}, we complete the proof.
\end{proof}

In the given context, we denote the average of an iterative sequence. Taking $\{x_k\}_{k=0}^{\infty}$ for example, its iterative average is given by
\[
\overline{x}_{N} = \frac{1}{N} \sum_{k=1}^{N}x_k.
\]
By utilizing~\Cref{lem: pdhg}, we can determine the rate that the average sequence $\{x_{k}\}_{k=0}^{\infty}$ converges in the weak sense. We state this result rigorously as follows.
\begin{theorem}
\label{thm: pdhg-weak}
Under the assumptions of~\Cref{lem: pdhg}, the average of the iterative sequence $\{\overline{x}_{N}\}_{N=1}^{\infty}$ converges weakly to the set $\mathrm{VI}$, denoted as~\eqref{eqn: vi}, at the rate characterized by
\begin{multline}
\label{eqn: pdhg-wek}
f(\overline{x}_{N}) - f(x) + g^{\star}(\overline{y}_{N}) - g^{\star}(y) + \big\langle F(\overline{x}_{N} - x), y \big\rangle  -  \big\langle Fx, \overline{y}_{N} - y \big\rangle \\ \leq \frac{\|x_0 - x\|^2 + \|y_0 - y\|^2 - 2s\big\langle F(x_0 - x), y_0 - y \big\rangle}{2sN}
\end{multline} 
for any initial $(x_0, y_0) \in \mathbb{R}^{d_1} \times \mathbb{R}^{d_2}$. 
\end{theorem}

When the point $(x,y)$ is identified as the saddle point $(x^{\star},y^{\star})$, the Lyapunov function~\eqref{eqn: pdhg-lyapunov} is reformulated as
\begin{equation}
\label{eqn: pdhg-saddle}
\mathcal{E}(k) = \frac{1}{2s} \|x_k - x^{\star}\|^2 + \frac{1}{2s} \|y_k - y^{\star}\|^2 - \big\langle F(x_k - x^{\star}), y_k - y^{\star} \big\rangle.
\end{equation}
Furthermore, the iterative difference~\eqref{eqn: pdhg-lem1} can be encapsulated by the ensuing inequality as
\begin{align}
\mathcal{E}(k+1) - &\mathcal{E}(k) \nonumber \\
                   \leq & -  \big\langle F(x_{k+1} - x^{\star}), y^{\star} \big\rangle  + \big\langle Fx^{\star} , y_{k+1} - y^{\star} \big\rangle +  f(x^{\star}) - f(x_{k+1}) + g^{\star}(y^{\star}) - g^{\star}(y_{k+1})  \nonumber \\
                           & - \underbrace{ \left(   \frac{1}{2s} \|x_{k+1} - x_{k}\|^2 +  \frac{1}{2s}  \|y_{k+1} - y_{k}\|^2 -  \big\langle F(x_{k+1} - x_{k}), y_{k+1} - y_{k} \big\rangle \right)}_{\textbf{NE}}. \label{eqn: pdhg-inq-saddle}
\end{align}
Building on~\Cref{thm: saddle-sn}, it becomes evident that the iterative difference, as captured by the inequality~\eqref{eqn: pdhg-inq-saddle}, can be negatively bounded by the numerical error resulting from the implicit scheme, while the corresponding time derivative, delineated by the inequality~\eqref{eqn: high-ode-inq-saddle}, remains only non-positive for the continuous scenario. This observation allows us to establish the convergence rates for the discrete~\texttt{PDHG} algorithm, articulated rigorously as follows.
\begin{theorem}
\label{thm: pdhg-ergodic}
Under the assumptions of~\Cref{lem: pdhg}, the iterative sequence $\{(x_{k}, y_{k})\}_{k=0}^{\infty}$ converges to $(x^{\star}, y^{\star})$ at the following rates characterized by
\begin{subequations}
\begin{empheq}[left=\empheqlbrace]{align}
          \frac{1}{N+1} \sum_{k=0}^{N} \big(  \|x_{k+1} - x_k\|^2 +  &\|y_{k+1} - y_{k}\|^2 - 2s\big\langle F(x_{k+1} - x_{k}), y_{k+1} - y_{k} \big\rangle \big) \nonumber \\ 
          & \leq \frac{\|x_0 - x^{\star}\|^2 + \|y_0 - y^{\star}\|^2 - 2s\big\langle F(x_0 - x^{\star}), y_0 - y^{\star} \big\rangle}{N+1}, \label{eqn: pdhg-ergodic-aver} \\
           \min_{0\leq k\leq N}\big(  \|x_{k+1} - x_k\|^2 +  \|&y_{k+1} - y_{k}\|^2 - 2s\big\langle F(x_{k+1} - x_{k}), y_{k+1} - y_{k} \big\rangle \big) \nonumber \\ 
          & \leq \frac{\|x_0 - x^{\star}\|^2 + \|y_0 - y^{\star}\|^2 - 2s\big\langle F(x_0 - x^{\star}), y_0 - y^{\star} \big\rangle}{N+1},   \label{eqn: pdhg-ergodic-min}
\end{empheq}
\end{subequations} 
for any initial $(x_0, y_0) \in \mathbb{R}^{d_1} \times \mathbb{R}^{d_2}$. Particularly, when the step size satisfies $s \in \big(0, \|F\|^{-1}\big)$, the iterative sequence $\{(x_{k}, y_{k})\}_{k=0}^{\infty}$ converges strongly with the rates as
\begin{subequations}
\begin{empheq}[left=\empheqlbrace]{align}
          \frac{1}{N+1} \sum_{k=0}^{N} \big(  \|x_{k+1} - x_k\|^2 +  \|y_{k+1} - y_{k}\|^2  \big) &\leq \frac{\left(1+s\|F\| \right)\left(\|x_0 - x^{\star}\|^2 + \|y_0 - y^{\star}\|^2 \right)}{ \left(1 - s\|F\| \right) (N+1)}, \label{eqn: pdhg-ergodic1-aver} \\
          \min_{0\leq k\leq N}\big(  \|x_{k+1} - x_k\|^2 +  \|y_{k+1} - y_{k}\|^2 \big)  & \leq \frac{ \left(1+s\|F\| \right)\left(\|x_0 - x^{\star}\|^2 + \|y_0 - y^{\star}\|^2 \right) }{\left(1-s \|F\| \right) (N+1)}.   \label{eqn: pdhg-ergodic1-min}
\end{empheq}
\end{subequations}                  
for any initial $(x_0, y_0) \in \mathbb{R}^{d_1} \times \mathbb{R}^{d_2}$.
\end{theorem}

Similarly, if we further assume that the objective function satisfies $f \in \mathcal{S}_{\mu,L}^{1}(\mathbb{R}^d)$, the convex inequality~\eqref{eqn: subgrad2-pdhg-descent} can be refined. This results in a tighter inequality as
\begin{equation}
\label{eqn: pdhg-strong-inq}
 f(x) - f(x_{k+1}) \geq \left\langle \frac{x_{k+1} - x_{k}}{s} - F^{T}(y_{k+1} - y_{k}) + F^{\top}y_{k+1}, x_{k+1} - x\right\rangle + \frac{\mu}{2} \|  x_{k+1} - x \|^2. 
\end{equation}
Following the process that leads to~\Cref{thm: pdhg-ergodic}, we can derive that the iterative difference also satisfies the following inequality as
\begin{align}
\mathcal{E}(k+1) - \mathcal{E}(k) \leq &  -  \frac{\mu}{2} \|  x_{k+1} - x^{\star} \|^2 \nonumber \\
                           & - \underbrace{ \left(   \frac{1}{2s} \|x_{k+1} - x_{k}\|^2 +  \frac{1}{2s}  \|y_{k+1} - y_{k}\|^2 -  \big\langle F(x_{k+1} - x_{k}), y_{k+1} - y_{k} \big\rangle \right)}_{\textbf{NE}}. \label{eqn: pdhg-inq-saddle-strong}
\end{align}
By incorporating the convexity of the $\ell_2$-norm square, we can determine the convergence rate of the average of the iterative sequence as stated in the following theorem.
\begin{theorem}
\label{thm: pdhg-strong}
Under the assumptions of~\Cref{lem: pdhg}, if we further assume the objective function $f \in \mathcal{S}_{\mu,L}^{1}(\mathbb{R}^d)$, the average of the iterative sequence $\{\overline{x}_{N}\}_{N=0}^{\infty}$ converges strongly to the saddle $(x^{\star}, y^{\star})$ at the rate characterized by
\begin{equation}
\label{eqn: pdhg-strong}
\|\overline{x}_{N} - x^{\star}\|^2 \leq  \frac{\|x_0 - x^{\star}\|^2 + \|y_0 - y^{\star}\|^2 - 2s\big\langle F(x_0 - x^{\star}), y_0 - y^{\star} \big\rangle}{\mu s N}
\end{equation}
for any initial $(x_0, y_0) \in \mathbb{R}^{d_1} \times \mathbb{R}^{d_2}$.
\end{theorem}

\section{Monotonicity}
\label{sec: monotonicity}

In this section, we aim to explore the monotonic behavior of the numerical error, as shown in~\eqref{eqn: pdhg-inq}. The numerical error is given by
\begin{equation}
\label{eqn: numerical-error}
\mathbf{NE} =  \frac{1}{2s} \|x_{k+1} - x_{k}\|^2 +  \frac{1}{2s}  \|y_{k+1} - y_{k}\|^2 -  \big\langle F(x_{k+1} - x_{k}), y_{k+1} - y_{k} \big\rangle.
\end{equation}
Our goal is to investigate how this error evolves as the iteration $k$ progresses. A monotonic decrease in the numerical error ($\mathbf{NE}$) with successive iterations indicates that the convergence rate $O(1/N)$, as referenced in~\Cref{thm: pdhg-ergodic}, can be improved for the last iterate. This potential enhancement is of notable significance, especially when compared with the convergence rates achieved through averaging and minimization strategies.


To begin our investigation,  it is instructive to consider the continuous counterpart of the numerical error~\eqref{eqn: numerical-error}, conceptualized as a Lyapunov function denoted by
\begin{equation}
\label{eqn: lyapunov-mono-ode}
\mathcal{E}(t) = \frac{1}{2s} \|\dot{X}\|^2 + \frac{1}{2s} \|\dot{Y}\|^2 - \big\langle F\dot{X}, \dot{Y} \big\rangle.
\end{equation}
Employing the classical Lyapunov analysis to take the time derivative of~\eqref{eqn: lyapunov-mono-ode}, we have
\begin{equation}
\label{eqn: lyapunov-mono-ode1}
\frac{d\mathcal{E}}{dt} = \frac{1}{s} \big\langle \ddot{X} - sF^{\top}\ddot{Y}, \dot{X} \big\rangle + \frac{1}{s} \big\langle \ddot{Y} - sF\ddot{X}, \dot{Y} \big\rangle. 
\end{equation}
Assuming that both functions, $f \in \mathcal{F}^{2}(\mathbb{R}^{d_1}) \cap \mathcal{R}(\mathbb{R}^{d_1}) $ and $g \in \mathcal{F}^{2}(\mathbb{R}^{d_1})$, are sufficiently smooth, we take the time derivative of each ODE in the high-resolution system,~\eqref{eqn: high-descent} and~\eqref{eqn: high-ascent}, to obtain 
\begin{subequations}
\begin{empheq}[left=\empheqlbrace]{align}
  & \ddot{X} - s F^{\top} \ddot{Y} =-  F^{T}\dot{Y} - \nabla^2 f(X)\dot{X},                    \label{eqn: high-descent-derivative} \\
  & \ddot{Y} - s F\ddot{X} = F\dot{X} - \nabla^2 g^{\star}(Y)\dot{Y}.                             \label{eqn: high-ascent-derivative}
\end{empheq}
\end{subequations}
Substituting these expressions,~\eqref{eqn: high-descent-derivative} and~\eqref{eqn: high-ascent-derivative}, into the derivitive inequality~\eqref{eqn: lyapunov-mono-ode1}, we conclude that the time derivative satisfies
\begin{equation}
\label{eqn: lyapunov-mono-ode2}
\frac{d\mathcal{E}}{dt} = -\big\langle \nabla^2 f(X)\dot{X}, \dot{X} \big\rangle - \big\langle \nabla^2 g^{\star}(Y)\dot{Y}, \dot{Y} \big\rangle \leq 0,
\end{equation}
which indicates that the Lyapunov function~\eqref{eqn: lyapunov-mono-ode1} decreases montonically.

Building on the insights from the system of high-resolution ODEs,~\eqref{eqn: high-descent} and~\eqref{eqn: high-ascent}, we now proceed to analyze the discrete~\texttt{PDHG} algorithm,~\eqref{eqn: pdhg-descent} ---~\eqref{eqn: pdhg-ascent}.  As indicated in~\eqref{eqn: lyapunov-mono-ode}, the numerical error $\mathbf{NE}$ is conceptualized within this framework as the discrete Lyapunov function, expressed as
\begin{equation}
\label{eqn: lyapunov-mono-pdhg}
\mathcal{E}(k) = \frac{1}{2s} \left\|x_{k+1} - x_{k}\right\|^2 + \frac{1}{2s} \left\|y_{k+1} - y_{k}\right\|^2- \big\langle F(x_{k+1} - x_{k}), y_{k+1} - y_{k} \big\rangle. 
\end{equation}
To analyze the iterative behavior of the Lyapunov function~\eqref{eqn: lyapunov-mono-pdhg}, we observe  its change across iterations as
\begin{align}
             \mathcal{E}(k&+1)  -  \mathcal{E}(k) \nonumber\\
                           \leq & \left\langle \frac{x_{k+2} - 2x_{k+1} + x_{k}}{s} - F^{\top}\left[y_{k+2} - 2y_{k+1} + y_{k})\right], x_{k+2} - x_{k+1} \right\rangle \nonumber \\
                                  & + \left\langle \frac{y_{k+2} - 2y_{k+1} + y_{k}}{s} - F\left[x_{k+2} - 2x_{k+1} + x_{k}) \right], y_{k+2} - y_{k+1} \right\rangle \nonumber \\
                               = &  \left\langle \frac{x_{k+2} - 2x_{k+1} + x_{k}}{s} - F^{\top}\left[y_{k+2} - 2y_{k+1} + y_{k})\right] + F^{\top}(y_{k+2}- y_{k+1}), x_{k+2} - x_{k+1} \right\rangle \nonumber \\
                                   & + \left\langle \frac{y_{k+2} - 2y_{k+1} + y_{k}}{s} - F\left[x_{k+2} - 2x_{k+1} + x_{k}) \right] - F(x_{k+2} - x_{k+1}) , y_{k+2} - y_{k+1} \right\rangle. \label{eqn: eqn: lyapunov-mono-pdhg-diff}
\end{align}
By invoking~\Cref{defn: subgradient} and the two raltions,~\eqref{eqn: subgrad1-pdhg-descent} and~\eqref{eqn: subgrad1-pdhg-ascent}, we can derive the following two inequalities as
\begin{subequations}
\begin{empheq}[left=\empheqlbrace]{align}
  & \left\langle \frac{x_{k+2} - 2x_{k+1} + x_{k}}{s} - F^{\top}\left[y_{k+2} - 2y_{k+1} + y_{k})\right] + F^{\top}(y_{k+2}- y_{k+1}), x_{k+2} - x_{k+1} \right\rangle \geq 0,                    \label{eqn: pdhg-descent-diff} \\
  & \left\langle \frac{y_{k+2} - 2y_{k+1} + y_{k}}{s} - F\left[x_{k+2} - 2x_{k+1} + x_{k}) \right] - F(x_{k+2} - x_{k+1}) , y_{k+2} - y_{k+1} \right\rangle \geq 0.                           \label{eqn: pdhg-ascent-diff}
\end{empheq}
\end{subequations}
Subsequently, by substituting these inequalities,~\eqref{eqn: pdhg-descent-diff} and~\eqref{eqn: pdhg-ascent-diff}, into the inequality of the iterative difference~\eqref{eqn: eqn: lyapunov-mono-pdhg-diff}, we establish that the Lyapunov function~\eqref{eqn: lyapunov-mono-pdhg} satisfies
\[
\mathcal{E}(k+1)  -  \mathcal{E}(k) \leq 0,
\]
which confirms that the numerical error $\mathbf{NE}$~\eqref{eqn: numerical-error} diminishes monotonically. This result allows us to refine~\Cref{thm: pdhg-ergodic}, enhancing the articulation of the convergence rate, especially as relevant for the last iterate. The rigorous statement is given as follows.

\begin{theorem}
\label{thm: pdhg-mono}
Under the assumptions of~\Cref{lem: pdhg}, the iterative sequence $\{(x_{k}, y_{k})\}_{k=0}^{\infty}$ converges to $(x^{\star}, y^{\star})$ at the following rates characterized by
\begin{multline}
\label{eqn: pdhg-mono}
            \|x_{k+1} - x_k\|^2 +  \|y_{k+1} - y_{k}\|^2 - 2s\big\langle F(x_{k+1} - x_{k}), y_{k+1} - y_{k} \big\rangle \\   \leq \frac{\|x_0 - x^{\star}\|^2 + \|y_0 - y^{\star}\|^2 - 2s\big\langle F(x_0 - x^{\star}), y_0 - y^{\star} \big\rangle}{N+1}, 
\end{multline}
for any initial $(x_0, y_0) \in \mathbb{R}^{d_1} \times \mathbb{R}^{d_2}$. Particularly, when the step size satisfies $s \in \big(0, 1/ \|F\|\big)$, the iterative sequence $\{(x_{k}, y_{k})\}_{k=0}^{\infty}$ converges strongly with the rates as
\begin{equation}
\label{eqn: pdhg-mono1}
          \|x_{k+1} - x_k\|^2 +  \|y_{k+1} - y_{k}\|^2  \leq \frac{\left(1+s\|F\| \right)\left(\|x_0 - x^{\star}\|^2 + \|y_0 - y^{\star}\|^2 \right)}{ \left(1 - s\|F\| \right) (N+1)},   
\end{equation}
for any initial $(x_0, y_0) \in \mathbb{R}^{d_1} \times \mathbb{R}^{d_2}$.
\end{theorem}

%% file: 05_extension.tex
\section{Extension: The general form of PDHG}
\label{sec: general-pdhg}

From~\Cref{sec: ode} to~\Cref{sec: monotonicity}, we have elaborated on the framework of high-resolution ODEs to decipher the convergence behavior of the discrete~\texttt{PDHG} algorithm,~\eqref{eqn: pdhg-descent} ---~\eqref{eqn: pdhg-ascent}. It is clear that its effectiveness hinges on the coupled corrections, $x$-correction and $y$-correction. Particularly, we make use of Lyapunov functions,~\eqref{eqn: pdhg-lyapunov} and~\eqref{eqn: lyapunov-mono-pdhg}, to derive its convergence rates. In this section, we extend this framework to include the general form of~\texttt{PDHG}, which is expressed as
\begin{subequations}
\begin{empheq}[left=\empheqlbrace]{align}
  & x_{k+1} = \argmin_{x \in \mathbb{R}^{d_1}} \left\{ f(x) + \langle Fx, y_k\rangle + \frac{1}{2\tau} \|x - x_k \|^2 \right\},                                            \label{eqn: pdhg-ex-descent} \\
  & \overline{x}_{k+1} = x_{k+1} + (x_{k+1} - x_k),                                                                                                                                                              \label{eqn: pdhg-ex-special} \\
  & y_{k+1} = \argmax_{y \in \mathbb{R}^{d_2}} \left\{ -g^{\star}(y) + \langle F\overline{x}_{k+1}, y\rangle - \frac{1}{2\sigma} \|y - y_k \|^2 \right\}.  \label{eqn: pdhg-ex-ascent}
\end{empheq}
\end{subequations}
where the two parameters $\tau$ and $\sigma$ satisfy the condition $0 < \tau\sigma\|F\|^{2} \leq 1$.  In other words, the general form,~\eqref{eqn: pdhg-ex-descent} ---~\eqref{eqn: pdhg-ex-ascent}, introduces two distinct parameters, $\tau$ and $\sigma$, for the proximal updates,  diverging from~\texttt{PDHG}~\eqref{eqn: pdhg-descent} ---~\eqref{eqn: pdhg-ascent} that utilizes  a single step size $s$.

Let the step size be $s = \sqrt{\tau\sigma}$. Following the same processes outlined in~\Cref{subsec: derivation-ode}, we can derive a system of high-resolution ODEs for the general form of~\texttt{PDHG},~\eqref{eqn: pdhg-ex-descent} ---~\eqref{eqn: pdhg-ex-ascent}, as 
\begin{subequations}
\begin{empheq}[left=\empheqlbrace]{align}
  & \alpha\dot{X} - s F^{\top} \dot{Y} = - F^{T}Y - \nabla f(X)                  \label{eqn: high-ex-descent} \\
  & \beta\dot{Y} - s F\dot{X} =  FX - \nabla g^{\star}(Y)                           \label{eqn: high-ex-ascent}
\end{empheq}
\end{subequations}
where the two parameters  $\alpha$ and $\beta$ satisfy the condition $\alpha \beta = 1$. This condition indicates that the term $( \alpha\dot{X},  \beta\dot{Y} )^{\top}$ represents the principal term, being of at least order $O(1)$, and thus play a dominant role in the dynamics. It is worth noting that the $O(s)$ term, $ (- s F^{\top} \dot{Y}, - s F\dot{X})^{\top}$, serves as a coupled perturbation within the high-resolution system,~\eqref{eqn: high-ex-descent} and~\eqref{eqn: high-ex-ascent}.

Accordingly, the Lyapunov function~\eqref{eqn: lyapunov-ode-saddle} is adapted to the following form as
\begin{equation}
\label{eqn: lyapunov-ex-ode}
\mathcal{E}(t) = \frac{1}{2\tau} \|X - x^{\star}\|^2 + \frac{1}{2\sigma} \|Y - y^{\star}\|^2 - \big\langle F(X - x^{\star}), Y - y^{\star} \big\rangle. 
\end{equation}
With the help of~\Cref{thm: subgrad-unique} and~\Cref{thm: saddle-sn}, it is straightforward for us to show that its time derivative is non-positive. Furthermore, the discrete Lyapunov function, is initially formulated in~\eqref{eqn: pdhg-saddle}, is modified to 
\begin{equation}
\label{eqn: pdhg-ex-lyapunov}
\mathcal{E}(k) = \frac{1}{2\tau} \|x_k - x^{\star}\|^2 + \frac{1}{2\sigma} \|y_k - y^{\star}\|^2 - \big\langle F(x_k - x^{\star}), y_k - y^{\star} \big\rangle. 
\end{equation}
Applying a similar approach that used in establishing~\Cref{thm: pdhg-ergodic}, we deduce that the iterative difference in iterations is bounded by the numerical error as 
\begin{equation}
\label{eqn: numerical-error-ex}
\mathbf{NE} =  \frac{1}{2\tau} \|x_{k+1} - x_{k}\|^2 +  \frac{1}{2\sigma}  \|y_{k+1} - y_{k}\|^2 -  \big\langle F(x_{k+1} - x_{k}), y_{k+1} - y_{k} \big\rangle.
\end{equation}
Therefore, we can enhance~\Cref{thm: pdhg-ergodic} as the following statement.
\begin{theorem}
\label{thm: pdhg-ergodic-general}
Let $f \in \mathcal{F}^{0}(\mathbb{R}^{d_1}) \cap \mathcal{R}(\mathbb{R}^{d_1}) $ and $g \in \mathcal{F}^{0}(\mathbb{R}^{d_1})$. If the two parameters $\tau$ and $\sigma$ satisfy  $\tau\sigma \in \big(0,  \|F\|^{-2}\big]$,  the iterative sequence $\{(x_{k}, y_{k})\}_{k=0}^{\infty}$ converges to $(x^{\star}, y^{\star})$ at the following rates characterized by
\begin{subequations}
\begin{empheq}[left=\empheqlbrace]{align}
          \frac{1}{N+1} \sum_{k=0}^{N} \big(  \sigma \|x_{k+1} - x_k\|^2 +  &\tau\|y_{k+1} - y_{k}\|^2 - 2\tau\sigma\big\langle F(x_{k+1} - x_{k}), y_{k+1} - y_{k} \big\rangle \big) \nonumber \\ 
          & \leq \frac{\sigma \|x_0 - x^{\star}\|^2 + \tau \|y_0 - y^{\star}\|^2 - 2 \tau\sigma \big\langle F(x_0 - x^{\star}), y_0 - y^{\star} \big\rangle}{N+1}, \label{eqn: pdhg-ergodic-aver-gen} \\
           \min_{0\leq k\leq N}\big(  \sigma \|x_{k+1} - x_k\|^2 +  \tau &\|y_{k+1} - y_{k}\|^2 - 2 \tau \sigma \big\langle F(x_{k+1} - x_{k}), y_{k+1} - y_{k} \big\rangle \big) \nonumber \\ 
          & \leq \frac{ \sigma \|x_0 - x^{\star}\|^2 + \tau \|y_0 - y^{\star}\|^2 - 2 \tau \sigma\big\langle F(x_0 - x^{\star}), y_0 - y^{\star} \big\rangle}{N+1},   \label{eqn: pdhg-ergodic-min-gen}
\end{empheq}
\end{subequations} 
for any initial $(x_0, y_0) \in \mathbb{R}^{d_1} \times \mathbb{R}^{d_2}$. Particularly,  when the two parameters $\tau$ and $\sigma$ satisfy $\tau\sigma \in \big(0,  \|F\|^{-2}\big)$, the iterative sequence $\{(x_{k}, y_{k})\}_{k=0}^{\infty}$ converges strongly with the rates as
\begin{subequations}
\begin{empheq}[left=\empheqlbrace]{align}
          \frac{1}{N+1} \sum_{k=0}^{N} \big(  \|x_{k+1} - x_k\|^2 +  \|y_{k+1} - y_{k}\|^2  \big) &\leq \frac{\left(1+\sqrt{\tau\sigma}\|F\| \right)\left(\|x_0 - x^{\star}\|^2 + \|y_0 - y^{\star}\|^2 \right)}{ \left(1 - \sqrt{\tau\sigma}\|F\| \right) (N+1)}, \label{eqn: pdhg-ergodic1-aver-gen} \\
          \min_{0\leq k\leq N}\big(  \|x_{k+1} - x_k\|^2 +  \|y_{k+1} - y_{k}\|^2 \big)  & \leq \frac{ \left(1+\sqrt{\tau\sigma}\|F\| \right)\left(\|x_0 - x^{\star}\|^2 + \|y_0 - y^{\star}\|^2 \right) }{\left(1-\sqrt{\tau\sigma} \|F\| \right) (N+1)}.   \label{eqn: pdhg-ergodic1-min-gen}
\end{empheq}
\end{subequations}                  
for any initial $(x_0, y_0) \in \mathbb{R}^{d_1} \times \mathbb{R}^{d_2}$.
\end{theorem}

Finally, we confirm the monotonic decrease of the numerical error~\eqref{eqn: numerical-error-ex}, which aligns with the process as described in~\Cref{sec: monotonicity}. This starts by considering the continuous system of high-resolution ODEs,~\eqref{eqn: high-ex-descent} and~\eqref{eqn: high-ex-ascent}. By taking the continuous counterpart of the numerical error~\eqref{eqn: numerical-error-ex} as a Lyapunov function, denoted as
\begin{equation}
\label{eqn: lyapunov-ex-mono-ode}
\mathcal{E}(t) = \frac{1}{2\tau} \|\dot{X}\|^2 + \frac{1}{2\sigma} \|\dot{Y}\|^2 - \big\langle F\dot{X}, \dot{Y} \big\rangle,
\end{equation}
it is straightforward for us to deduce that the time derivative of $\mathcal{E}(t)$ in~\eqref{eqn: lyapunov-ex-mono-ode} does not exceed zero. This conclusion is achieved by taking the time derivative of the system of high-resolution ODEs,~\eqref{eqn: high-ex-descent} and~\eqref{eqn: high-ex-ascent}. Proceeding with a similar approach for the discrete Lyapunov function
\begin{equation}
\label{eqn: lyapunov-ex-mono-pdhg}
\mathcal{E}(k) = \frac{1}{2\tau} \left\|x_{k+1} - x_{k}\right\|^2 + \frac{1}{2\sigma} \left\|y_{k+1} - y_{k}\right\|^2 - \big\langle F(x_{k+1} - x_{k}), y_{k+1} - y_{k} \big\rangle. 
\end{equation}
Through some elementary analysis, we can discern that the numerical error indeed decreases monotonically. This insight allows us to refine~\Cref{thm: pdhg-mono} accordingly.
\begin{theorem}
\label{thm: pdhg-mono-general}
Let $f \in \mathcal{F}^{0}(\mathbb{R}^{d_1}) \cap \mathcal{R}(\mathbb{R}^{d_1}) $ and $g \in \mathcal{F}^{0}(\mathbb{R}^{d_1})$. If the two parameters $\tau$ and $\sigma$ satisfy $\tau\sigma \in \big(0,  \|F\|^{-2}\big]$,  the iterative sequence $\{(x_{k}, y_{k})\}_{k=0}^{\infty}$ converges to $(x^{\star}, y^{\star})$ at the following rates characterized by
\begin{multline}
\label{eqn: pdhg-mono-gen}
            \sigma \|x_{k+1} - x_k\|^2 + \tau \|y_{k+1} - y_{k}\|^2 - 2\tau\sigma \big\langle F(x_{k+1} - x_{k}), y_{k+1} - y_{k} \big\rangle \\   \leq \frac{ \sigma \|x_0 - x^{\star}\|^2 +  \tau \|y_0 - y^{\star}\|^2 - 2\tau\sigma \big\langle F(x_0 - x^{\star}), y_0 - y^{\star} \big\rangle}{N+1}, 
\end{multline}
for any initial $(x_0, y_0) \in \mathbb{R}^{d_1} \times \mathbb{R}^{d_2}$. Particularly,  when the two parameters $\tau$ and $\sigma$ satisfy $\tau\sigma \in \big(0,  \|F\|^{-2}\big)$, the iterative sequence $\{(x_{k}, y_{k})\}_{k=0}^{\infty}$ converges strongly with the rates as
\begin{equation}
\label{eqn: pdhg-mono1-gen}
         \sigma \|x_{k+1} - x_k\|^2 + \tau  \|y_{k+1} - y_{k}\|^2  \leq \frac{\left(1+\sqrt{\tau\sigma}\|F\| \right)\left(\sigma\|x_0 - x^{\star}\|^2 + \tau \|y_0 - y^{\star}\|^2 \right)}{ \left(1 - \sqrt{\tau\sigma}\|F\| \right) (N+1)},   
\end{equation}
for any initial $(x_0, y_0) \in \mathbb{R}^{d_1} \times \mathbb{R}^{d_2}$.
\end{theorem}


%% file: 06_conclu.tex
\section{Conclusion and discussion}
\label{sec: conclu}

In this study, we utilize the dimensional analysis, a method previously employed in~\citep{shi2022understanding, shi2023learning, shi2021hyperparameters, li2024understanding}, to derive a system of high-resolution ODEs for the discrete~\texttt{PDHG} algorithm. This system effectively captures a key feature of~\texttt{PDHG}, the coupled $x$-correction and $y$-correction, distinguishing it from the proximal Arrow-Hurwicz algorithm. The small but essential perturbation ensures that~\texttt{PDHG} consistently converges, bypassing the periodic behavior observed in the proximal Arrow-Hurwicz algorithm. Technically, we utilize Lyapunov analysis to investigate the convergence behavior as it shifts from the continuous high-resolution system to the discrete~\texttt{PDHG} algorithm. This analysis identifies the numerical errors, resulting from the implicit scheme, as the crucial factor that leads to the convergence rate and monotonicity in the discrete~\texttt{PDHG} algorithm, echoing a significant observation also made in the context of the \texttt{ADMM} algorithm as identified in~\citep{li2024understanding}. In addition, we further discover that if the objective function $f$ is strongly convex, the iterative average of~\texttt{PDHG} is enhanced and converges strongly.

Understanding and analyzing the discrete~\texttt{PDHG} algorithm via the system of high-resolution ODEs opens up several exciting avenues for further investigation. In comparison with the system of high-resolution ODEs of~\texttt{ADMM} as outlined in~\citep[(1.16a) --- (1.16c)]{li2024understanding}, the high-resolution system of~\texttt{PDHG},~\eqref{eqn: high-descent} and~\eqref{eqn: high-ascent}, features a dynamics structure that is both simpler and more intuitive. It would be compelling to explore the convergence of~\texttt{PDHG} across various norms and rates. Furthermore, given the second-order gradient system that simulates a heavy ball accelerating down the valley to speed up the gradient flow in a manner akin to the motion of a droplet, there is considerable potential to devise an algorithm enhancing the convergence rates of~\texttt{PDHG} to reach the lower bound $O(1/N^{2})$ for the convex optimization challenges, as highlighted in~\citep{nemirovski1983problem}. Turning our attention back to the objective function~\eqref{eqn: game-rep}, and excluding the two single-variable functions, the main convex-concave function blending two variables is quadratic. Building on the intuitive dynamics of the high-resolution system of~\texttt{PDHG},~\eqref{eqn: high-descent} and~\eqref{eqn: high-ascent}, it would also be appealing to extend this framework to non-quadratic functions, particularly to explore the convergence behavior of the~\texttt{PDHG} algorithm on non-quadratic convex-concave objective functions. Additionally, a very intriguing research direction is to investigate the traditional minimax optimization algorithms, such as~\textit{optimal gradient descent ascent} (\texttt{OGDA}) algorithm~\citep{popov1980modification} and the Extragradient method~\citep{korpelevich1976extragradient} and their relationship with~\texttt{PDHG} through high-resolution ODEs.
